\documentclass{amsart}
\usepackage{amssymb,latexsym,enumerate,amscd}
\newtheorem{Thm}{Theorem}[section]
\newtheorem{Prop}[Thm]{Proposition}
\newtheorem{Lem}[Thm]{Lemma}
\newtheorem{Cor}[Thm]{Corollary}

\newtheorem{Conj}[Thm]{Conjecture}

\theoremstyle{remark}
\newtheorem{Rem}[Thm]{Remark}
\newtheorem{Cond}[Thm]{Condition}
\newtheorem{Ass}[Thm]{Assertion}
\newtheorem*{Cav}{Caveat}
\theoremstyle{definition}
\newtheorem{Not}[Thm]{Notation}
\newtheorem*{Not*}{Notation}

\numberwithin{equation}{section}
\DeclareMathOperator{\im}{im}
\DeclareMathOperator{\coker}{coker}
\DeclareMathOperator{\tr}{tr}
\DeclareMathOperator{\Tr}{Tr}
\DeclareMathOperator{\hatTr}{\widehat{Tr}}
\DeclareMathOperator{\id}{id}
\DeclareMathOperator{\GL}{GL}
\DeclareMathOperator{\Wh}{Wh}
\DeclareMathOperator{\Iso}{Iso}
\DeclareMathOperator{\Hom}{Hom}
\DeclareMathOperator{\Orthog}{Orthog}

\DeclareMathOperator{\Mat}{Mat}
\DeclareMathOperator{\Units}{Units}
\newcommand{\dummy}{\,.\,}
\begin{document}

\date{%
Sat Sep 14 19:11:23 EDT 2002}

\title[Whitehead groups and the Bass conjecture]
{Whitehead groups and the Bass conjecture}

\author[F. T. Farrell]{F. Thomas Farrell}
\address{Department of Mathematical Sciences \\
State University of New York at Binghamton \\
Binghamton \\
NY 13902-6000 \\
USA}
\email{farrell@math.binghamton.edu}
\urladdr{http://www.math.binghamtom.edu/farrell/}
\thanks{The first author was supported in part by the National
Science Foundation}

\author[P. A. Linnell]{Peter A. Linnell}
\address{Department of Mathematics \\
Virginia Tech \\
Blacksburg \\
VA 24061-0123 \\
USA}
\email{linnell@math.vt.edu}
\urladdr{http://www.math.vt.edu/people/linnell/}

\begin{abstract}
This paper will be concerned with proving that
certain Whitehead groups of
torsion-free elementary amenable groups are torsion groups and
related results, and then applying these results to the Bass
conjecture.  In particular we shall establish the strong Bass
conjecture for an arbitrary elementary amenable group.
\end{abstract}



\maketitle

\section{Introduction} \label{Sintroduction}

Let $k$ be a field, let $\Gamma$ be a group,
and let $\Units^*(k\Gamma)$ denote the subgroup of
$\Units(k\Gamma)$ consisting of those units in $k\Gamma$ of the form
$a\gamma$ with $a \in k \setminus 0$ and $\gamma \in \Gamma$.
Then we shall use the notation
$\Wh^k(\Gamma)$ for the quotient of $K_1(k\Gamma)$ by the
image of $\Units^*(k\Gamma)$ under the natural homomorphism
$\Units^*(k\Gamma) \to K_1(k\Gamma)$.  
All rings will have a unity element 1, and a ring of prime
characteristic will mean a ring such that $p1 = 0$ for some prime
$p$.  Let $\mathcal {T}$ indicate the
class of torsion abelian groups.  Recall that the class of
elementary amenable groups is the smallest class of groups which
contains $\mathbb {Z}$ and all finite groups, and is closed under
taking extensions and directed unions.  Our first result is

\begin{Thm} \label{Twhitehead}
Let $k$ be a field of prime characteristic and let
$\Gamma$ be a torsion-free elementary amenable group.  Then
$\Wh^k(\Gamma) \in \mathcal {T}$.
\end{Thm}

Let $\tilde {K}_0(R)$ denote the reduced projective class
group of the ring $R$,
that is $K_0(R)/\langle [R] \rangle$.  Then in the situation of
Theorem \ref{Twhitehead}, it now follows from Bass's
contracted functor argument \cite[Chapter XII, \S 7]{Bass68}
that $\tilde{K}_0(k\Gamma) \in \mathcal {T}$
and $K_i(k\Gamma) \in \mathcal
{T}$ for all $i < 0$, where $k$ and $\Gamma$ are as in Theorem
\ref{Twhitehead}.  However we shall prove the following
more general result.

\begin{Thm} \label{TK_0}
Let $k$ be a field of prime characteristic and let $\Gamma$ be an
elementary amenable group.  Then $K_0(k\Gamma) \otimes \mathbb {Q}$
is generated by the images of $K_0(kG)$ as $G$ runs over the finite
subgroups of $\Gamma$.
\end{Thm}

It is likely that Theorems \ref{Twhitehead} and
\ref{TK_0} are also true for fields $k$ of characteristic 0.  However
our proofs depend on Lemma \ref{Lkey}, which shows that for rings of
prime characteristic, certain nil groups are torsion groups.  This is
certainly false for arbitrary rings, so our proof does not cover the
case when $k$ has characteristic 0.

Given a prime power $q$,
we shall let $k_q$ indicate the field with $q$ elements.
In the case
$k = k_2$, we shall give a second proof of Theorem
\ref{TK_0} in Section \ref{Sappendix}.
This alternative proof will yield other results;
for example we obtain the following variation of \cite[Theorem
1.1]{FarrellLinnell}.
\begin{Thm} \label{Tlinear}
Let $G$ be a torsion-free virtually solvable subgroup of
$\GL_n(\mathbb {C})$.  Then $\Wh^{k_2}(G) = 0$.
\end{Thm}
\noindent
We shall also prove (a 2-group is a group in which
every element has order a power of 2)
\begin{Prop} \label{Pwreath}
Let $G$ be an abelian  2-group.  Then $\tilde{K}_0(
k_2[G \wr \mathbb {Z}]) = 0$.
\end{Prop}
For the special case $G=C_2$ in Proposition \ref{Pwreath}, the group
is often called the lamplighter group.  Many interesting properties
of this group are established in \cite{GrigorchukZuk01}.

We shall apply Theorem \ref{TK_0} to obtain results on Bass's strong
conjecture \cite[4.5]{Bass79}.  For convenience, we will restate the
conjecture here (see Section \ref{Spreliminary} for the
explanation of some standard notation used below).
\begin{Conj}[The strong Bass conjecture] \label{CBass}
Let $k$ be an integral domain, let $\Gamma$ be a group, let $g \in
\Gamma$, and let $P$ be a finitely generated projective
$k\Gamma$-module.  Suppose $o(g)$ is not invertible in
$k$.  Then $r_P(g) = 0$.
\end{Conj}
Some cases for which
Conjecture \ref{CBass} is known to be true are
\begin{enumerate}
\item $k = \mathbb
{C}$ and $G$ a linear group \cite[Proposition 6.2]{Bass79},
\item
$k = \mathbb {Z}$ and $o(g) < \infty$ \cite[Lemma 4.1]{Linnell83}, 
\item
$k = \mathbb {Q}$ and $G$ has cohomological dimension
over $\mathbb {Q}$ at most two \cite[Theorem 3.3]{Eckmann86}.
\end{enumerate}
We shall prove
\begin{Thm} \label{TBass}
The strong Bass conjecture is true for elementary amenable groups.
More precisely,
let $\Gamma$ be an elementary amenable group, let $k$ be an
integral domain, and let $P$ be a finitely generated projective
$k\Gamma$-module.
If $g \in \Gamma$ and $o(g)$
is not invertible in $k$, then $r_P(g) = 0$.
\end{Thm}

While revising this paper, we have learned that Berrick, Chatterji
and Mislin \cite{BerrickChatterjiMislin}
have proved the strong Bass conjecture for amenable groups in
the case $k= \mathbb {C}$.  In fact they prove a rather stronger
version for the Banach space $\ell^1(G)$; on the other hand their
results do not apply to the case when $k$ has nonzero characteristic.
Their techniques are rather different from ours, and depend on
recent work of Vincent Lafforgue \cite{Lafforgue}.

One of many interesting results which Gerald Cliff proved in his
important paper \cite{Cliff80} was \cite[Theorem 1]{Cliff80},
that if $k$ is a field of nonzero characteristic
$p$ and $\Gamma$ is a polycyclic-by-finite group with the
property that all finite subgroups have $p$-power order, then
$k\Gamma$ has no nontrivial idempotents.
We will use Theorem \ref{TK_0} and
the techniques of Section \ref{SBass} to
extend Cliff's result to elementary amenable groups.
\begin{Thm} \label{Tidemp}
Let $p$ be a prime, let $k$ be an integral
domain of characteristic
$p$, and let $\Gamma$ be an elementary amenable group.  Suppose
every finite subgroup of $\Gamma$ has $p$-power order.
Then $k\Gamma$ has no nontrivial idempotents.
\end{Thm}

This work was carried out while we were at the
Sonderforschungsbereich in M\"unster.  We would like to thank
Wolfgang L\"uck for organizing our visits to M\"unster, and the
Sonderforschungsbereich for financial support.

\section{Preliminary results} \label{Spreliminary}

\begin{Not*}
All modules will be right modules and mappings will be written on the
left.  For each positive integer $n$, we shall use the notation
$C_n$ for the cyclic group of order $n$, $\Mat_n(R)$ for the $n
\times n$ matrices over the ring $R$, $\Mat(R)$ for
$\bigcup_{n=1}^{\infty} \Mat_n(R)$, and
$\mathcal {A}$ for the class of finitely
generated virtually abelian groups.
Also $G\wr A$ will denote the restricted wreath product
of the groups $G$ and $A$; thus $G\wr \mathbb {Z}$ will have
base group $\bigoplus_{i \in \mathbb {Z}} G$ and quotient group
$\mathbb {Z}$.  We shall let $o(g)$ denote the order of the element
$g \in G$.  In the case $o(g) = \infty$, we shall adopt the
convention that $o(g)$ is \emph{not} invertible in any ring.
A $\mathcal {T}$-exact sequence will mean a sequence which is exact
modulo torsion abelian groups; in other words every element of
the homology group at each stage of the sequence has finite order.
Similarly a $\mathcal {T}$-epimorphism means a
homomorphism which is onto modulo torsion.
Suppose $\alpha$ is an automorphism of the ring $R$.  Then
$R_{\alpha}[t]$ will denote the twisted polynomial ring over $R$, and
$R_{\alpha}[t,t^{-1}]$ will denote the twisted Laurent polynomial
ring over $R$.  Following
\cite[\S 2]{FarrellHsiang70}, we define $\mathcal {C}(R, \alpha)$
to be the category whose objects are pairs $(P, \phi)$ where $P$
is a finitely generated projective $R$-module and $\phi$ is an
$\alpha$-linear nilpotent endomorphism of $P$, and whose morphisms $g
\colon (P_1, \phi_1) \to (P_2, \phi_2)$ are $R$-linear
homomorphisms $g \colon P_1 \to P_2$ with $g\phi_1 = \phi_2 g$.  Then
we shall let $C(R, \alpha) = K_0(\mathcal {C}(R, \alpha))/\langle
[(R,0)] \rangle$, and $\tilde {C}(R, \alpha)$ denote the subgroup
of $C(R, \alpha)$ generated by elements of the form $[(R^n, \phi)]$.
We remark that $\tilde {C}(R, \alpha)$ is isomorphic to the subgroup
of $K_1(R_{\alpha}[t])$ generated by elements which are
represented by matrices in $\Mat(R_{\alpha}[t])$
of the form $I + Nt$, where $I$ is the
identity matrix and $N \in \Mat(R)$ such that $Nt$ is nilpotent;
this can be seen from \cite[proof of Theorem 13 and
Proposition 20]{FarrellHsiang70}.
If $\beta$ is an automorphism of the group
$\Gamma$, then we shall let $\beta$ also indicate the automorphism
of $k\Gamma$ induced by $\beta$.
\end{Not*}

We need a standard description
(cf.\ \cite[\S 3]{KrophollerLinnell88})
of elementary amenable groups in order
to carry out induction arguments.  If $\mathcal {X}$ and $\mathcal
{Y}$ are classes of groups, then $G \in \mathcal {X}
\mathcal {Y}$ will mean that
$G$ has a normal subgroup $H$ such that $H \in \mathcal {X}$ and $G/H
\in \mathcal {Y}$, and $G \in L\mathcal {X}$ will mean that every
finitely generated subgroup of $G$ is contained in a $\mathcal
{X}$-group (if $\mathcal {X}$ is closed under taking subgroups, this
is equivalent to saying that every finitely generated subgroup of $G$
is an $\mathcal {X}$-group).  For each ordinal $\alpha$, we define
$\mathcal {X}_{\alpha}$ inductively as follows.
\begin{verse}
$\mathcal {X}_0$ is the class of finite groups;

$\mathcal {X}_{\alpha} = (L \mathcal {X}_{\alpha -1})\mathcal {A}$ if
$\alpha$ is a successor ordinal;

$\mathcal {X}_{\alpha} = \bigcup_{\beta < \alpha} \mathcal
{X}_{\beta}$ if $\alpha$ is a limit ordinal.
\end{verse}
Then the proof of \cite[Lemma 3.1]{KrophollerLinnell88}
(see also \cite[Lemma 4.9]{Linnell93}) yields the following.
\begin{Lem} \label{LEA}
\begin{enumerate} [\normalfont (i)]
\item
The class of elementary amenable groups is $\bigcup _{\alpha \ge 0}
\mathcal {X}_{\alpha}$;
\item
Each $\mathcal {X}_{\alpha}$ is
closed under taking subgroups;
\item
If $H \lhd G$ with $G/H$ finite and $H \in
\mathcal {X}_{\alpha}$ or $L\mathcal {X}_{\alpha}$,
then $G \in \mathcal {X}_{\alpha}$ or $L\mathcal {X}_{\alpha}$
respectively.
\end{enumerate}
\end{Lem}

\begin{Lem} \label{Lkey}
Let $p$ be a prime, let $R$ be a ring such that $p1 = 0$, and let
$\alpha$ be an automorphism of $R$.  Then
$\tilde {C} (R, \alpha) $ is an abelian $p$-group.
\end{Lem}
\begin{proof}
Let $t$ be an indeterminant and let $N \in \Mat(R)$
such that $(Nt)^n = 0$ for some positive integer
$n$.  Then $I + Nt$ represents an element of $K_1( 
R_{\alpha}[t])$ and we need to prove that this element has
finite order.  Now for $M \in \Mat (R_{\alpha}[t])$,
the binomial formula shows that
$(I+M)^p = I+M^p$ because $pM = 0$ and $p$ divides
$\binom{p}{i}$ when $0 < i < p$.  This equation by repetition yields
\[
(I+M)^{p^n} = I + M^{p^n}.
\]
Substituting into this equation $M = Nt$ and observing that
$(Nt)^{p^n} = 0$ because $p^n \ge n$, we obtain
$(I + Nt)^{p^n} = I$.  Hence $\tilde {C} (R, \alpha)$ is
generated by elements of $p$-power order and the result follows.
\end{proof}

\begin{Lem} \label{LC}
Let $R$ be a ring of prime characteristic and let $\alpha$ be an
automorphism of $R$.
Then there is a natural $\mathcal {T}$-exact sequence
\[
K_1(R) \longrightarrow K_1( R_{\alpha}[t,t^{-1}]) \longrightarrow
K_0(R) \overset{1 - \alpha_*}{\longrightarrow} K_0(R).
\]
\end{Lem}
\begin{proof}
This follows immediately from \cite[Theorem 19c]{FarrellHsiang70} and
Lemma \ref{Lkey}.
\end{proof}

\begin{Rem} \label{Rbass}
The constant group homomorphism $\mathbb {Z}^n \to 1$ induces a ring
homomorphism
\[
R[\mathbb {Z}^n] \to R[1] = R
\]
which splits the inclusion ring homomorphism $R\to R[\mathbb {Z}^n]$.
Hence $K_i(R[\mathbb {Z}^n])$ is naturally a direct sum of $K_i(R)$
and the complementary summand $\hat{K}_i(R[\mathbb {Z}^n])$ which is
the kernel of the homomorphism induced by the above splitting.  The
Bass-Heller-Swan formula, which is the fundamental theorem of
algebraic $K$-theory, asserts that
\[
\hat{K}_1(R[\mathbb {Z}]) \cong K_0(R) \oplus C(R,\id) \oplus
C(R,\id).
\]
In his contracted functor theory Bass used this
formula to define the lower
algebraic $K$-groups $K_{-n}(R)$ for $n > 0$, so that they are direct
summands of $\hat{K}_0(R[\mathbb {Z}^n])$; cf.\ \cite[Chapter XII, \S
7]{Bass68}.
\end{Rem}

\begin{Cor} \label{Cnil1}
Let $k$ be a field of prime characteristic and let $\Gamma = \pi
\rtimes \mathbb {Z}$ be a group such that $\Wh^k(\pi \times \mathbb
{Z}^n) \in \mathcal {T}$ for all $n \ge 0$.  Then $\Wh^k(\Gamma
\times \mathbb {Z}^n) \in \mathcal {T}$ for all $n \ge 0$.
\end{Cor}
\begin{proof}
Note that \cite[Theorem 21d]{FarrellHsiang70} remains true with $k$
in place of $\mathbb {Z}$ and $\Wh^k$ in place of $\Wh$.  So applying
this result, Lemma \ref{Lkey} and the Bass-Heller-Swan formula,
we obtain the exact sequence
\begin{multline*}
0 = \Wh^k(\pi \times \mathbb {Z}^n) \otimes
\mathbb {Q} \longrightarrow
\Wh^k(\Gamma \times \mathbb {Z}^n) \otimes \mathbb {Q}
\longrightarrow \tilde
{K}_0(k[\pi \times \mathbb {Z}^n]) \otimes \mathbb {Q} \\
\subseteq \Wh^k(\pi \times
\mathbb {Z}^{n+1}) \otimes \mathbb {Q} = 0.
\end{multline*}
This proves $\Wh^k(\Gamma \times \mathbb {Z}^n) \in \mathcal {T}$ for
all $n \ge 0$, as required.
\end{proof}

\begin{Cor} \label{Cnil2}
Let $R$ be a ring of prime characteristic and let
$\alpha$ be an automorphism of $R$.  Suppose
the natural map $K_0(R) \to K_0(R[s,s^{-1}])$ is a $\mathcal
{T}$-epimorphism.  Then the natural map $K_0(R) \to K_0(R_{\alpha}
[t,t^{-1}])$ is also a $\mathcal {T}$-epimorphism.
\end{Cor}
\begin{proof}
Consider the following commutative diagram
\[
\hspace*{-3ex}
\begin{CD}
K_0(R) @>>> K_0(R_{\alpha}[t,t^{-1}]) \\
@VV
\parbox{1ex}{$\uparrow$\\$\hspace*{.35ex}
\overset{\textstyle\shortmid}{\shortmid}$}
V  @VV
\parbox{1ex}{$\uparrow$\\$\hspace*{.35ex}
\overset{\textstyle\shortmid}{\shortmid}$}
V  \\
K_1(R[s,s^{-1}]) @>>> K_1(R_{\alpha}[s,s^{-1},t,t^{-1}])
@>>> K_0(R[s,s^{-1}])
@>1-\alpha_*>>K_0(R[s,s^{-1}])\\
@AAA @AAA @A\mathcal{T}\text{-epi}AA
@AA
\parbox{1ex}
{$\hspace*{.35ex}\underset{\textstyle\shortmid}{\shortmid}$\\
$\downarrow$}
A\\
K_1(R) @>>> K_1(R_{\alpha}[t,t^{-1}])
@>>> K_0(R)
@>1-\alpha_*>> K_0(R)
\end{CD}
\]
The squares in this diagram all commute and the 3 pairs of vertical
arrows going in the opposite directions are splittings.  The two
horizontal sequences are $\mathcal {T}$-exact by Lemma \ref{LC}.
Also in each of
the first two columns, the composite of the two up vertical
arrows is 0.  A
simple diagram chase using the fact that the marked up vertical arrow
is a $\mathcal {T}$-epi yields the result.
\end{proof}

\begin{Cor} \label{Cnil3}
Let $R$ be a ring of prime characteristic,
let $\alpha$ be an automorphism of $R$, and suppose that the
natural map $K_0(R) \to K_0(R[\mathbb {Z}^n])$ is a
$\mathcal {T}$-epimorphism for all $n$.  
Then $K_i(R) \in \mathcal {T}$ and $K_i(R_{\alpha}[t,t^{-1}][\mathbb
{Z}^n]) \in \mathcal {T}$ for all $i<0$, and $K_0(R)
\to K_0(R_{\alpha}[t,t^{-1}][\mathbb {Z}^n])$ is a $\mathcal
{T}$-epimorphism, for all $n$.
\end{Cor}
\begin{proof}
We use Bass's contracted functor theory to deduce this; in
particular, we use Remark \ref{Rbass}.  By our assumption
$\hat{K}_0(R[\mathbb {Z}^n]) \in \mathcal{T}$ and hence $K_i(R) \in
\mathcal{T}$ for all $i<0$ because $K_i(R) \subseteq
\hat{K}_0(R[\mathbb {Z}^n])$ where $n = -i$.  Since
\[
K_0(R[\mathbb {Z}^n]) \longrightarrow K_0(R[\mathbb {Z}^{n+1}]) =
K_0(R[\mathbb {Z}^n][s,s^{-1}])
\]
is clearly a $\mathcal{T}$-epimorphism for all $n\ge 0$, $f_1 \colon
K_0(R[\mathbb {Z}^n]) \to K_0(R_{\alpha}[t,t^{-1}][\mathbb {Z}^n])$
is also a $\mathcal{T}$-epimorphism by Corollary \ref{Cnil2}.
Consequently so is $f_2 \colon K_0(R) \to
K_0(R_{\alpha}[t,t^{-1}][\mathbb {Z}^n])$ since it is the composite
of the two $\mathcal{T}$-epimorphisms $K_0(R) \to K_0(R[\mathbb
{Z}^n])$ and $f_1$.  Setting $n = m + (-i)$, where $m \ge 0$ and $0 >
i$ are given, we see that
\[
f_3 \colon K_0(R_{\alpha}[t,t^{-1}][\mathbb {Z}^m]) \longrightarrow
K_0((R_{\alpha}[t,t^{-1}][\mathbb {Z}^m])[\mathbb {Z}^{-i}]),
\]
is also a $\mathcal{T}$-epimorphism since $f_2 = f_3 \circ f_4$ where
$f_4 \colon K_0(R) \to K_0(R_{\alpha}[t,t^{-1}][\mathbb {Z}^m])$.
Therefore $\hat{K}_0(S[\mathbb {Z}^{-i}]) \in \mathcal{T}$ where $S =
R_{\alpha}[t,t^{-1}][\mathbb {Z}^m]$.  But $K_i(S)$ is a direct
summand of $\hat{K}_0(S[\mathbb {Z}^i])$, consequently
$K_i(R_{\alpha}[t,t^{-1}][\mathbb {Z}^m]) \in \mathcal{T}$.
\end{proof}

\begin{Cor} \label{Cnil4}
Let $k$ be a field of prime characteristic, let
$\Gamma = \pi \rtimes \mathbb {Z}$
be a group such that for all $n \ge 0$,
$K_0(k[\pi \times \mathbb {Z}^n]) \otimes \mathbb {Q}$ is generated
by the images of $K_0(kG)$ as $G$ varies over the finite subgroups of
$\pi$.  Then for all $i<0$ and $n \ge 0$, it follows that
$K_i(k\pi) \in \mathcal {T}$ and
$K_i(k[\Gamma \times \mathbb {Z}^n]) \in \mathcal {T}$, and
$K_0(k[\Gamma \times \mathbb {Z}^n]) \otimes \mathbb {Q}$ is
generated by the images of $K_0(kG)$ as $G$ varies over the finite
subgroups of $\Gamma$.
\end{Cor}
\begin{proof}
Let $\alpha$ denote the automorphism of $\pi$ determined by the
conjugation action of $\mathbb {Z}$ on $\pi$, and let $R = k\pi$.
Note that
\begin{align*}
R[\mathbb {Z}^n] &= k[\pi \times \mathbb {Z}^n] \\
\text{and} \quad
R_{\alpha}[t,t^{-1}][\mathbb {Z}^n] &= k[\Gamma \times \mathbb
{Z}^n].
\end{align*}
Thus the natural map $K_0(R) \to K_0(R[\mathbb {Z}^n])$ is
$\mathcal {T}$-surjective, and we can now obtain the result
from Corollary \ref{Cnil3}.
\end{proof}

\begin{Lem} \label{Lfinitegroup}
Let $n \ge 0$, let
$G$ be a finite group, and let $k$ be a field.
Then $K_i(kG) = 0$ for all $i < 0$ and
the inclusion $G \hookrightarrow G \times \mathbb {Z}^n$ induces
an isomorphism $K_0(kG) \to K_0(k[G \times \mathbb {Z}^n])$.
\end{Lem}
\begin{proof}
Let $J$ denote the Jacobson radical of $kG$.  Since $kG$ is an
Artinian ring, $J$ is nilpotent \cite[Corollary 1, p.~39]{Jacobson64}
and $kG/J$ is a semisimple Artinian ring \cite[\S III.3]{Jacobson64},
and in particular $kG/J$ is a regular ring.
Then we have a commutative diagram
\[
\begin{CD}
K_0(kG) @>>> K_0(k[G \times \mathbb {Z}^n]) \\
@V\cong VV @VV\cong V \\
K_0(kG/J) @>\cong>>
K_0(k[G \times \mathbb {Z}^n]/J[\mathbb {Z}^n]).
\end{CD}
\]
The two vertical arrows are isomorphisms by
\cite[Lemma II.2.2]{Weibel} and the
bottom horizontal line is an isomorphism by the Bass-Heller-Swan
formula, hence the top arrow is also an isomorphism as required and
the second statement is proven.

A consequence of the second statement is that $\hat{K}_0(S[\mathbb
{Z}^{-i}]) = 0$ where $S = kG$ and $i < 0$.  Hence $K_i(kG) = 0$
since it is a direct summand of $\hat{K}_0(kG[\mathbb {Z}^{-i}])$.
\end{proof}

\begin{Lem} \label{Lliegroup}
Let $L \in \mathcal {A}$.  Then $L$ is
isomorphic to a discrete cocompact subgroup of a
virtually connected Lie group.
\end{Lem}
\begin{proof}
There clearly exists a non-negative integer $n$ such that $L$ is an
extension of a finite group $F$ by the free abelian group $\mathbb
{Z}^n$, where $\mathbb {Z}^n$ denotes the integral lattice points in
the additive group $\mathbb {R}^n$.  This extension determines an
action of $F$ on $\mathbb {Z}^n$ and hence on $\mathbb {R}^n$, and a
cohomology class $\theta \in H^2(F,\mathbb {Z}^n)$.  Let $\theta' \in
H^2(F,\mathbb {R}^n)$ be the image of $\theta$ and let $G$ be the
extension of $F$ by $\mathbb {R}^n$ determined by $\theta'$.  Then
$G$ is a virtually connected Lie group containing $L$ as a discrete
cocompact subgroup.  (Note that $\theta' = 0$ and hence
$G = \mathbb {R}^n \rtimes F$.)
\end{proof}

The subclass of $\mathcal{A}$ consisting of the virtually cyclic
groups is of particular importance to us.  There is fortunately the
following quite useful structure theorem for this subclass due to
Scott and Wall \cite{ScottWall79};
cf.\ \cite[Lemma 2.5]{FarrellJones95} for another proof.

\begin{Prop} \label{Pvircyclic}
A virtually cyclic group $\Gamma$ contains a finite normal subgroup
$F$ such that $\Gamma/F$ is either trivial, infinite cyclic, of
infinite dihedral.
\end{Prop}

\section{Whitehead groups of elementary amenable groups}

\begin{Thm} \label{Tfiber}
Let $k$ be a field of prime characteristic, and let
$\pi \lhd \Gamma$ be groups such that $\Gamma/\pi$ is a
crystallographic group.
Suppose $\Wh^k(\tilde {\pi} \times \mathbb {Z}^n) \in \mathcal {T}$
for all non-negative integers $n$ whenever
$\tilde {\pi}/\pi $ is a finite subgroup of $\Gamma/\pi$.  Then
$\Wh^k(\Gamma) \in \mathcal {T}$.
\end{Thm}

\begin{proof}
Our proof of Theorem \ref{Tfiber} follows the pattern established in
\cite{FarrellHsiang81, Quinn85}.  In these papers it was shown
that $\Wh^R(\Gamma) = 0$ for $\Gamma$ a torsion-free virtually
poly-$\mathbb {Z}$ group and $R$ any subring of $\mathbb {Q}$.
The case $R = \mathbb {Z}$ was done in \cite{FarrellHsiang81} and
the general case in \cite{Quinn85}.  Quinn had to develop important
new geometric algebra concepts to do the general case, and these
concepts are crucial in our proof of Theorem \ref{Tfiber}.

\begin{Not} \label{Nfiber1}
Let $L$ be a crystallographic group.  Then $L$ is isomorphic to a
discrete cocompact subgroup of the group of all rigid motions
\[
\Iso (\mathbb {E}^n) \cong \Orthog (n) \ltimes \mathbb {R}^n
\]
of some Euclidean space $\mathbb {E}^n$.  The number $n$ is the rank
of a torsion-free abelian subgroup of finite index in $L$, and is
called the dimension of $L$ or $\dim L$.  Also the image of $L$ in
$\Orthog (n)$ is a finite group $G$, called the
\emph{holonomy group} of $L$,
and its isomorphism class is determined by $L$.  The order of $G$ is
called the \emph{holonomy number} of $L$ and is denoted by $\#(L)$.
If $S$ is a subgroup of finite index in $L$, then $S$ is also
crystallographic.  Furthermore $\dim S = \dim L$ and the holonomy
group $G_1$ of $S$ is isomorphic to a subgroup of $G$, hence
\begin{itemize}
\item If $G$ is cyclic, then so is $G_1$.
\item $\# (S) \le \#(L)$.
\end{itemize}
\end{Not}

Let $L = \Gamma/\pi$ where $\Gamma$ and $\pi$ come from the
statement of Theorem \ref{Tfiber}, let $\phi \colon \Gamma
\twoheadrightarrow L$ denote the natural epimorphism,
and let $G$ denote the holonomy
group of $L$.  Frobenius induction relative to $G$ reduces the proof
of Theorem \ref{Tfiber} to the case $G$ is cyclic.  
This follows from \cite[Corollary 2.12]{Swan70} (see
\cite{FarrellHsiang81} for more details).
It is also recommended that the reader now glance at the summary of
Swan's ``Frobenius induction theory" given later in this paper in the
paragraph following \eqref{Eorder}; in particular, see the important
fact \eqref{Eswan3} mentioned there.
We proceed to prove Theorem \ref{Tfiber} by
simultaneous induction on $\dim (L)$ and $\#(L)$, where we always
assume that $G$ is cyclic.  Our explicit inductive assumption is that
whenever $L_0$ is a crystallographic group such that either
\begin{alignat*}{2}
\dim(L_0) & < \dim(L) &\quad& \text{or}\\
\dim(L_0) &= \dim(L) &\quad&\text{and } \#(L_0) < \#(L),
\end{alignat*}
then the theorem is true for $L_0$.  So primary induction on $\dim(L)$
and secondary induction on $\#(L)$.

To start the $n$th secondary induction, we need to show that Theorem
\ref{Tfiber} is true when $L= \mathbb {Z}^n$ (i.e.\ when $\#(L) =
1$).  Clearly we may assume
that $n > 0$.  In this case $\Gamma/\pi \cong \mathbb {Z}^n$, so
there exists a normal subgroup $\Gamma_0$ of $\Gamma$ containing
$\pi$ such that $\Gamma/\Gamma_0 \cong \mathbb {Z}$ and $\Gamma_0/\pi
\cong \mathbb {Z}^{n-1}$.  Then $\Gamma \cong \Gamma_0 \rtimes
\mathbb {Z}$, and by induction $\Wh^k(\Gamma_0 \times \mathbb {Z}^m)
\in \mathcal {T}$ for all $m \ge 0$.  We may now apply Corollary
\ref{Cnil1} with $\pi = \Gamma_0$.

\begin{Rem} \label{RZ}
Let $T$ be the abelian normal subgroup of the crystallographic group
$L$ consisting of all pure translations, so $T = L \cap \mathbb
{R}^n$ (see the notation above).  Then $L$ is an extension of $G$ by
$T$.  This extension determines by conjugation a representation of
$G$ on $T$, called the \emph{holonomy representation} of $L$.  It is
well known that $L$ maps epimorphically onto $\mathbb {Z}$
if and only if $T^G \ne 0$, where $T^G$ denotes the subgroup fixed by
$G$: cf.\ \cite[Lemma 1.4]{FarrellHsiang78} and \cite{Farkas75}.
This fact can be proven by applying the Lyndon-Hochschild-Serre
spectral sequence to the group extension $1 \to T \to L \to G \to 1$.
Its $E_2^{pq}$ term is $H^p(G;H^q(T,\mathbb {Z}))$ and it converges to
$H^{p+q}(L,\mathbb {Z})$.  Using that $E_2^{2,0}=H^2(G,\mathbb {Z})$
is finite, one sees easily that $E_2^{0,1} = (\Hom(L,\mathbb {Z}))^G$
vanishes if and only if $H^1(L,\mathbb {Z}) = \Hom(L,\mathbb {Z})$
vanishes.  But $(\Hom(T,\mathbb {Z}))^G$ vanishes if and only if
$T^G$ vanishes.
\end{Rem}

We now consider the general inductive step, which is divided into two
cases according to whether $T^G \ne 0$ or $T^G = 0$.
\bigskip

\noindent
Case 1. $T^G \ne 0$.
\newline
By Remark \ref{RZ}, we may write $L = \mathcal
{L} \rtimes \mathbb {Z}$ for some $\mathcal {L} \lhd L$.  Let
$\Gamma_0 = \phi^{-1}(\mathcal {L})$, so $\Gamma_0/\pi \cong \mathcal
{L}$.  Since an $\mathcal{A}$-group is
crystallographic if and only if its unique maximal finite normal
subgroup is 1, we see that $\mathcal {L}$ is a crystallographic group
with $\dim \mathcal {L} = \dim L -1$.
Therefore $\Wh^k(\Gamma_0 \times \mathbb {Z}^m) \in \mathcal {T}$
for all $m \ge 0$
because of our inductive assumption.  We may now apply Corollary
\ref{Cnil1} to complete the inductive step in Case 1.
\bigskip

\noindent
Case 2.  $T^G = 0$.
\newline
The geometric algebra developed by Quinn
\cite{Quinn79, Quinn82, Quinn85, Quinn85x} is used crucially
here, replacing the $h$-cobordisms used in \cite{FarrellHsiang81}.
A relatively simple example which concretely illustrates the
terminology used in the remainder of the proof of Theorem
\ref{Tfiber} is worked out in detain in Section \ref{Sappendix};
cf.\ Case 2 of the proof of Corollary \ref{CFIC} where $L$ is the
infinite dihedral group $C_2*C_2$.  It is recommended that the reader
keep this example in mind while perusing the rest of the proof of
Theorem \ref{Tfiber}.  He would also see \cite[Appendix]{Quinn82}
for details about stratified systems of fibrations.

Topology now enters into our proof.  Let $M$ be a connected smooth
manifold with $\pi_1(M) = \Gamma$ and denote its universal cover by
$\tilde {M}$.  Also identify $\Gamma$ with the group of all deck
transformations of $\tilde {M} \to M$.  Projection onto the first
factor of $E := \mathbb {R}^n \times_{\Gamma} \tilde {M}$ induces a
map
\[
q \colon E \longrightarrow \mathbb {R}^n/L,
\]
where $\Gamma$ acts on $\mathbb {R}^n$ via $\Gamma \twoheadrightarrow
L \subseteq \Iso (\mathbb {R}^n)$.
This map is a stratified system of fibrations on $\mathbb {R}^n/L$ in
the sense of \cite[Definition 8.2]{Quinn82}, whose strata are
determined in the standard way by the holonomy group action of $G$ on
the $n$-torus $\mathbb {R}^n/T$.  Let $s$ be any prime congruent to
1 mod $\#(L)$; recall that there is an infinitude of such primes by
Dirichlet's theorem.  For each such $s$, there is an endomorphism
$\psi$ of $L$ and a $\psi$-equivariant diffeomorphism $f \colon
\mathbb {R}^n \to \mathbb {R}^n$ (relative to $L \subseteq \Iso
(\mathbb {R}^n)$, so $f h = \psi (h) f$ for all $h \in L$)
such that
\begin{itemize}
\item
$| df(x) | = s | x | $ for each vector
$x$ tangent to $\mathbb {R}^n$;
\item $\psi(T) \subseteq T$;
\item $\psi$ induces $\id_G$ on $G=L/T$ and
\item $\psi |_T$ is multiplication by $s$.
\end{itemize}
This is a restating of \cite[Theorem 2.2]{FarrellHsiang81} which is
itself an immediate extension of the classical Epstein-Shub result
\cite{EpsteinShub68}.

\begin{Not} \label{Nfiber2}
We denote by $L_s$ and $T_s$ the finite quotient groups
$L/sT$ and $T/sT$ respectively.  Note that the exact sequence
\[
1 \longrightarrow T_s \longrightarrow L_s \longrightarrow G
\longrightarrow 1
\]
splits, and does so uniquely up to conjugacy since $(s, | G
| ) = 1$.
One splitting, which we shall denote by $G_s$, is given by the image
of $\psi (L)$ in $L_s$ under the projection $L \twoheadrightarrow
L_s$.  Let $\eta \colon \Gamma \to L_s$ denote
the composite epimorphism
\[
\Gamma \overset{\phi}{\twoheadrightarrow} L
\twoheadrightarrow L_s,
\]
and let $\Gamma_s \subseteq \Gamma$ indicate the inverse image
of $\psi (L)$ with respect to the epimorphism $\phi \colon \Gamma
\twoheadrightarrow L$, so $\Gamma_s/\pi = \psi (L)$.
Observe that $\Gamma_s = \eta^{-1}(G_s)$.
Furthermore, let $E_s$ denote
\[
\mathbb {R}^n \times_{\Gamma_s} \tilde{M},
\]
let $\hat{f} \colon \mathbb {R}^n/L \to \mathbb {R}^n/\psi(L)$ denote
the map induced by $f$,
and let $q_s \colon E_s \to \mathbb {R}^n/L$ be the composite of the
map $E_s \to \mathbb {R}^n/\psi(L)$ induced by projection onto the
first factor of $\mathbb {R}^n \times \tilde {M}$ with
the homeomorphism $\hat{f}^{-1}$.
\end{Not}

The map $q_s$ is also a stratified system of fibrations with the same
strata as $q$.  The following is an important observation.
\begin{Rem} \label{Rcurve}
Let $\alpha$ be a smooth curve in $E$ and let $\tilde {\alpha}$ be a
lift of $\alpha$ to the covering space $E_s$.  Then $| q_s \circ
\tilde {\alpha} | \le | q \circ \alpha |/s$, where
$| \dummy | $ denotes arc length (measured via $\mathbb {R}^n$).
\end{Rem}
\begin{proof}
Let $\hat{\alpha}$ be a lift of $\tilde{\alpha}$ to $\mathbb {R}^n
\times \tilde{M}$ which is the universal cover of both $E_s$ and $E$.
Then $\hat{\alpha} = (\hat{\alpha}_1,\hat{\alpha}_2)$ where
$\hat{\alpha}_1$ is a smooth curve in $\mathbb {R}^n$ and
$\hat{\alpha}_2$ is a smooth curve in $\tilde{M}$.  Then $|q_s
\circ \tilde{\alpha}|$ is by definition the arc length of $f^{-1}
\circ \alpha_1$ while $|q \circ \alpha|$ is the arc length of
$\alpha_1$.  Now the inequality (in fact equality) asserted in Remark
\ref{Rcurve} follows from the fact (noted above)
that $df$ stretches tangent vectors by a factor of $s$.
\end{proof}

Another important observation is the following.
\begin{Rem} \label{Rconjugate}
Let $S$ be a cyclic subgroup of $L_s$ such that $S$ projects onto
$G$ under the second map in the short exact sequence
\[
1 \longrightarrow T_s \longrightarrow L_s \longrightarrow G
\longrightarrow 1.
\]
Then $S \cap T_s = 1$; i.e., $S$ splits this sequence and is
consequently conjugate to $G_s$.
\end{Rem}
This observation is a consequence of our assumption that $T^G = 0$
together with \cite[Lemmas 1.2 and 1.4]{FarrellHsiang78}.

\begin{Rem} \label{Rsame}
So far most of the proof of Theorem \ref{Tfiber} can be repeated
verbatim for the proof of Theorem \ref{Tfiber0}
in the next section.  However at this
point the two proofs diverge somewhat.
\end{Rem}

We now recall some basic facts about Quinn assembly; cf.\
\cite[Appendix]{Quinn82} for more details.  These facts will also be
used in Sections \ref{S4} and
\ref{Sappendix}.  Let $\mathcal{S}$ be a homotopy
invariant (covariant) functor from the category of topological spaces
to $\Omega$-spectra.  Important examples of such functors are:
$X \mapsto \underline{K}(R\pi_1X)$,
$X \mapsto \underline{Wh}(\pi_1X)$,
$X \mapsto \underline{Wh}^k(\pi_1X)$,
where $R$ is a ring with 1.  Here $\underline{\Wh}(\pi)$ is the
cofiber of the standard map of spectra
\[
\mathbb{H}(B\pi;\underline{K}(\mathbb {Z}))
\longrightarrow
\underline{K}(\mathbb {Z}\pi)
\]
defined by Loday \cite{Loday76} and likewise $\underline{\Wh}^k(\pi)$
is the cofiber of
\[
\mathbb {H}(B\pi; \underline{K}(k))
\longrightarrow
\underline{K}(k\pi)
\]
defined also in \cite{Loday76}.
Let $\mathcal{M}$ denote the category of continuous
surjective maps; i.e.\ an
object in $\mathcal{M}$ is a continuous surjective map $p \colon E
\to B$ between topological spaces $E$ and $B$, while a morphism from
$p_1 \colon E_1 \to B_1$ to $p_2 \colon E_2 \to B_2$ is a pair of
continuous maps $f \colon E_1 \to E_2$, $g \colon B_1 \to B_2$ making
the following diagram a commutative square of maps:
\[
\begin{CD}
E_1 @>f>> E_2\\
@Vp_1VV @VVp_2V\\
B_1 @>>g> B_2.
\end{CD}
\]
Quinn \cite[Appendix]{Quinn82} constructed a functor from
$\mathcal{M}$ to the category of $\Omega$-spectra which associates to
the map $p$ the spectrum $\mathbb {H}(B;\mathcal{S}(p))$ in such a
way that $\mathbb{H}(B;\mathcal{S}(p)) = \mathcal{S}(E)$ in the
special case where $B$ is a single point *.  Furthermore the map of
spectra $\mathfrak{a} \colon \mathbb {H}(B;\mathcal{S}(p)) \to
\mathcal {S}(E)$ functorially associated to the commutative square
\[
\begin{CD}
E @>\id>> E\\
@VpVV @VVV\\
B @>>> *
\end{CD}
\]
is called the (Quinn) assembly map.  The map induced by
$\mathfrak{a}$ between the $i$th homotopy groups of these spectra is
also called the assembly map and is denoted by the same symbol;
i.e.\ $\mathfrak{a} \colon \mathbb{H}_i(B;\mathcal{S}(p))
\to \pi_i(\mathcal{S}(p))$.

Let $x$ be an arbitrary (but fixed) element of $\Wh^k(\Gamma)$.  We
need to show that $x$ has finite order to prove Theorem \ref{Tfiber}.
For this purpose, recall that Quinn showed that $x$ can be
represented by a geometric isomorphism (denoted $h$) of geometric
$k$-modules on the space $E$ \cite[\S 3.2]{Quinn85x}.  The transfer
of $h$ (denoted by $h_s$) to the finite sheeted cover $E_s \to
E$ clearly represents the transfer of $x$ to $\Wh^k(\Gamma_s)$
denoted by $x_s$.  But the radii $\epsilon_s$ of these geometric
isomorphisms $h_s$ measured via $q_s$ in $\mathbb {R}^n/L$ go to zero
as $s \to \infty$; i.e.
\begin{equation} \label{Eepsilon}
\lim_{s \to \infty}  \epsilon_s = 0.
\end{equation}
In fact $h_s$ represents an element $\bar{x}_s \in
\Wh^k(\mathbb{R}^n/L, q_s, \epsilon_s)$ and $\bar{x}_s$ maps to $x_s$
under the natural homomorphism
$\Wh^k(\mathbb {R}^n/L, q_s, \epsilon_s)
\to \Wh^k(\Gamma_s)$.  Since the strata in $\mathbb {R}^n/L$ of the
stratified systems $q_s$ are independent of $s$, Quinn's stability
theorem \cite[\S 4]{Quinn82} together with equation \eqref{Eepsilon}
and \cite[\S 3]{Quinn85} yield that $x_s$ is contained in the image
of $\mathbb {H}_1(\mathbb {R}^n/L; \underline{\Wh}^k(q_s))$ under the
Quinn assembly map provided $s$ is sufficiently large.  Now there is
an Atiyah-Hirzebruch-Quinn spectral sequence $E_{ij}^t$ converging to
$\mathbb {H}_{i+j}(\mathbb {R}^n/L; \underline{\Wh}^k(q_s))$ such that
\[
E^2_{ij} = H_i(\mathbb {R}^n/L; \Wh^k_j (q_s))
\]
\cite[Theorem 8.7]{Quinn82}.  Here
\[
\Wh_1^k(\dummy) = \Wh^k(\dummy), \quad
\Wh_0^k(\dummy) = \tilde{K}_0(\dummy), \quad
\Wh_j^k(\dummy) = K_j(\dummy)
\quad\text{if } j<0.
\]
The stalk of the coefficient sheaf $\Wh_j^k(q_s)$ over $y \in
\mathbb {R}^n/L$ is $\Wh_j^k \bigl(\pi_1 (q_s^{-1}(y))\bigr)$.
But $\pi$ is a subgroup of finite index in
$\pi_1 (q_s^{-1}(y))$.  Using the given hypotheses with $\tilde {\pi}
= \pi_1(q_s^{-1}(y))$ and Bass's contracted functor theory
\cite[Chapter XII, \S 7]{Bass68}, we see that
\[
\Wh_j^k \bigl(\pi_1(q_s^{-1}(y)) \bigr) \in \mathcal {T}
\quad\text{if } j\le 1.
\]
Since $\mathbb {R}^n/L$ is a finite polyhedron, we conclude that
$\mathbb {H}_1(\mathbb {R}^n/L; \underline{\Wh}^k(q_s))
\in \mathcal {T}$.  So this discussion yields that
\begin{equation} \label{Eorder}
x_s \text{ has finite order}
\end{equation}
provided $s$ is sufficiently large.  We now fix, for the remainder of
the proof, a sufficiently large prime $s$ such that
\begin{itemize}
\item $s \equiv 1 \mod | G |$ and
\item $x_s$ has finite order.
\end{itemize}
We proceed to apply Frobenius induction to $\Wh^k(\Gamma)$ relative
to the factor group $L_s$.  Let $\eta \colon \Gamma \to L_s$ denote
the composite epimorphism
\[
\Gamma \overset{\phi}{\twoheadrightarrow} L
\twoheadrightarrow L_s.
\]
Now $\Wh^k(\eta^{-1}(S))$ is a Frobenius
module over Swan's Frobenius functor
$G_0(S)$ as $S$ varies over the subgroups
of $L_s$ \cite{Swan70}.  Hence to show that $x$ has finite order, it
suffices to show, for each cyclic subgroup $S$ of $L_s$, that the
associated transfer of $x$ has finite order.  This criterion is a
consequence of \cite[Corollary 2.12]{Swan70}.  Furthermore, we need
only check this condition for one group in each conjugacy class of
cyclic subgroups of $L_s$.

For the reader's convenience, we recall the salient facts from Swan's
theory that are needed in this paper.  Swan associates to each finite
group $F$ a ring $G_0(F) = G_0(\mathbb {Z}F)$ with identity 1.  And
to each subgroup $S$ of $F$ he associates an additive group
homomorphism $G_0(S) \to G_0(F)$ called induction.  In particular
(after fixing $F$) $1_S$ denotes the image in $G_0(F)$ of $1 \in
G_0(S)$ under induction.  Let $\mathcal {C}$ be a collection of
cyclic subgroups of $F$ such that each conjugacy class is represented
exactly once in $\mathcal{C}$.  Then there exists a set of integers
$\{n(S) \mid S \in \mathcal{C}\}$ and a nonzero integer $n$ such that
\begin{equation} \label{Eswan1}
n = \sum_{S \in \mathcal{C}} n(S) 1_S.
\end{equation}
This is a consequence of \cite[Theorems 2.14, 2.19 and 4.1]{Swan70}.
Now let $\mathcal{S}(\pi)$ be any of the following functors from
groups $\pi$ to abelian groups:
\[
K_0(\mathbb {Z}\pi),\ K_0(k\pi),\ \tilde{K}_0(\mathbb {Z}\pi),
\ \tilde{K}_0(k\pi),\ K_1(\mathbb {Z}\pi),
\ K_1(k\pi), \Wh(\pi), \text{ or } \Wh^k(\pi).
\]
If $\psi \colon \Gamma \to F$ is a group epimorphism, then $A =
\mathcal{S}(\Gamma)$ is a (unitary) module over $G_0(F)$.  Define
abelian groups $A(S)$, for $S \in \mathcal{C}$, by $A(S) =
\mathcal{S}(\psi^{-1}(S))$.  Furthermore let $\sigma_S \colon A \to
A(S)$ and $\sigma^S \colon A(S) \to A$ denote the transfer and
induction group homomorphisms, respectively.  Then Swan proves the
following ``Frobenius reciprocity formula"
\begin{equation} \label{Eswan2}
(1_S)x = \sigma_S(\sigma^S(x))
\end{equation}
for all $x \in A$ and every $S \in \mathcal{C}$.  The following
important fact is an immediate consequence of formulae \eqref{Eswan1}
and \eqref{Eswan2}.
\begin{equation} \label{Eswan3}
\text{If } \sigma^S(x) \in \mathcal{T} \text{ for each }
S \in \mathcal{C}, \text{ then } x \in \mathcal{T}.
\end{equation}
And the following second important fact is also an easy consequence
of formulae \eqref{Eswan1} and \eqref{Eswan2}.
\begin{equation} \label{Eswan4}
\begin{minipage}{10.4cm}
Suppose for each $S \in \mathcal{C}$ that some nonzero (integral)
multiple of $\sigma^S(x)$ is the sum of elements induced from
$\mathcal{S}(H)$ as $H$ varies over the finite subgroups of
$\psi^{-1}(S)$.  Then some nonzero (integer) multiple of $x$ is also
a sum of elements induced from $\mathcal{S}(F)$ as $F$ varies over
the finite subgroups of $\Gamma$.
\end{minipage}
\end{equation}

If the natural projection $\sigma \colon L_s \twoheadrightarrow
G$ sends $S$ onto $G$, then $S$ is conjugate to $G_s$ by Remark
\ref{Rconjugate}.  But the transfer of $x$ associated to $G_s$ is
$x_s$, which has finite order because of \eqref{Eorder}.  We are
therefore left to examine the situation where $\sigma (S)$ is a
proper subgroup of $G$; i.e.\ we must look at the transfer of $x$ to
$\Wh^k(\eta^{-1}(S))$.  It clearly suffices to show that
$\Wh^k(\eta^{-1}(S)) \otimes \mathbb {Q} = 0$.  For this purpose,
consider the short exact sequence
\[
1 \longrightarrow \pi \longrightarrow \eta^{-1}(S)
\longrightarrow L(S) \longrightarrow 1
\]
where $L(S)$ denotes the inverse image of $S$ under the epimorphism
$L \twoheadrightarrow L_s$.
Observe that $L(S)$ is a crystallographic group  with
\begin{itemize}
\item $\dim L(S) = \dim L$ but
\item $\#(L(S)) < \#(L)$
\end{itemize}
because $\sigma (S) \ne G$.  Also notice that if $\hat{\pi}$ is any
subgroup of $\eta^{-1}(S)$ which contains $\pi$ with finite index,
then $\Wh^k(\hat{\pi} \times \mathbb {Z}^m) \otimes \mathbb {Q} = 0$
for every non-negative integer $m$ because $\eta^{-1}(S) \subseteq
\Gamma$ and we have assumed the same property for $\Gamma$.  We
therefore conclude from our inductive assumption that
$\Wh^k(\eta^{-1}(S)) \otimes \mathbb {Q} = 0$.
This completes the proof of Theorem \ref{Tfiber}.
\end{proof}

An easy consequence of Theorem \ref{Tfiber} is
\begin{Cor} \label{Cfibereasy}
Let $k$ be a field of prime characteristic, and let
$\pi \lhd \Gamma$ be groups such that $\Gamma/\pi$ is a
crystallographic group.
Suppose $\Wh^k(\tilde {\pi} \times \mathbb {Z}^n) \in \mathcal {T}$
for all non-negative integers $n$ whenever
$\tilde {\pi}/\pi $ is a finite subgroup of $\Gamma/\pi$.  Then
$\Wh^k(\Gamma \times \mathbb {Z}^n) \in \mathcal {T}$ for all $n \ge
0$.
\end{Cor}
\begin{proof}
Note that every finite subgroup of $\Gamma/\pi \times \mathbb {Z}^n$
is contained in $\Gamma/\pi$.
Since $\Gamma/\pi \times \mathbb {Z}^n$ is also a crystallographic
group, we can apply Theorem \ref{Tfiber} with $\Gamma = \Gamma \times
\mathbb {Z}^n$ and $\pi = \pi$.
\end{proof}

We now use Corollary \ref{Cfibereasy} to obtain the following result.
\begin{Cor} \label{Cfiber}
Let $k$ be a field of prime characteristic, and let
$\pi \lhd \Gamma$ be groups such that $\Gamma/\pi$ is an
elementary amenable group.
Suppose $\Wh^k(\tilde {\pi} \times \mathbb {Z}^n) \in \mathcal {T}$
for all non-negative integers $n$ whenever
$\tilde {\pi}/\pi $ is a finite subgroup of $\Gamma/\pi$.  Then
$\Wh^k(\Gamma \times \mathbb {Z}^n) \in \mathcal {T}$ for all $n \ge
0$.
\end{Cor}
Note that Theorem \ref{Twhitehead}
immediately follows from using Corollary
\ref{Cfiber} in the case $\pi = 1$.  To see this, since 1 is the only
finite subgroup of $\Gamma/1$, we need to show
that $\Wh^k(\mathbb {Z}^n) \in \mathcal {T}$ for all non-negative
integers $n$.  However we have in fact
$\Wh^k(\mathbb {Z}^n) = 1$ by the Bass-Heller-Swan
theorem, so Theorem \ref{Twhitehead} is proven.

\begin{proof}[Proof of Corollary \ref{Cfiber}]
We shall use the description of elementary amenable groups as
described in Lemma \ref{LEA}.
Let $\alpha$ be the least ordinal such that $\Gamma/\pi \in \mathcal
{X}_{\alpha}$.  Now $\alpha$ cannot be a limit ordinal, and the
result is clearly true if $\Gamma/\pi \in \mathcal {X}_0$ because
then $\Gamma/\pi$ is finite.  Therefore we may assume that $\alpha$
is a successor ordinal, and by transfinite induction that the result
is true for groups in $\mathcal {X}_{\alpha -1}$.
Thus $\Gamma$ has a normal subgroup $\pi_1$ containing
$\pi$ such that $\Gamma/\pi_1
\in \mathcal {A}$ and $\pi_1/\pi \in L\mathcal {X}_{\alpha - 1}$.
Now any $\mathcal {A}$-group maps onto a crystallographic group with
finite kernel; in other words there is a finite normal subgroup
$\pi_2/\pi_1$ of $\Gamma/\pi_1$ such that $\Gamma/\pi_2$ is
crystallographic.  Let
$\tilde {\pi}/\pi_2$ be any finite subgroup of $\Gamma/\pi_2$.
Then Lemma \ref{LEA} shows that $\tilde{\pi}/\pi
\in L\mathcal {X}_{\alpha -1}$.
Since Whitehead groups commute with direct limits, we see that
Corollary \ref{Cfiber} is true with $\tilde {\pi}$ in place of
$\Gamma$, so the result follows from Corollary \ref{Cfibereasy}.
\end{proof}

\section{$K_0$ of elementary amenable groups} \label{S4}

\begin{Thm} \label{Tfiber0}
Let $k$ be a field of prime characteristic and let
$\pi \lhd \Gamma$ be groups such that $\Gamma/\pi$
is a crystallographic group.  Suppose that
$K_0(k[\tilde {\pi} \times \mathbb {Z}^m]) \otimes \mathbb {Q}$ is
generated by the images of $K_0(kG) \otimes \mathbb {Q}$ as $G$
varies over the finite subgroups of $\tilde {\pi}$ where $\tilde
{\pi}/ \pi$ is a finite subgroup of $\Gamma / \pi$,
for all $m \ge 0$.  Then for every $m \ge 0$,
\begin{enumerate} [\normalfont (i)] \label{Tfiber0i}
\item $K_i(k[\Gamma \times \mathbb {Z}^m]) \otimes \mathbb {Q}
= 0$ for all $i < 0$.
\item \label{Tfiber0ii}
$K_0(k[\Gamma \times \mathbb {Z}^m])
\otimes \mathbb {Q}$ is generated by the images
of $K_0(kG) \otimes \mathbb {Q}$ as $G$
varies over the finite subgroups of $\Gamma$.
\end{enumerate}
\end{Thm}
\begin{proof}
Note that we need only establish assertion (ii), since Corollary
\ref{Cnil4} shows that (i) is a consequence of (ii).
Our proof of (ii) will follow that of Theorem \ref{Tfiber} up to
Remark \ref{Rsame} with minor modifications,
so we will keep Notation \ref{Nfiber1}.
We need to check that the $n$th secondary induction starts; i.e.\
Theorem \ref{Tfiber0} is true in the case $L = \mathbb {Z}^n$ (i.e.\
$\#(L) = 1$).  This is established by applying
Corollary \ref{Cnil4} with $\pi =
\Gamma_0$ (where $\pi \lhd \Gamma_0 \lhd \Gamma$ and $\Gamma/\Gamma_0
\cong \mathbb {Z}$).

We also keep Remark \ref{RZ} and the setup of Case 1.  However for
the last step in Case 1, we apply Corollary \ref{Cnil4} instead of
Corollary \ref{Cnil1}.

For Case 2, we keep Notation \ref{Nfiber2} and Remarks \ref{Rcurve}
and \ref{Rconjugate}.  We replace the
part of the proof of Theorem \ref{Tfiber} after Remark \ref{Rsame}
with the following new material.

We proceed to apply Frobenius induction relative
to the factor group $L_s$.
Let $x$ be an arbitrary but fixed element of $K_0(k[\Gamma \times
\mathbb {Z}^m])$.  For each cyclic subgroup $S$ of $L_s$, let $x(S)$
denote the transfer of $x$ to $K_0(k[\eta^{-1}(S) \times \mathbb
{Z}^n])$.  As $S$ varies over the subgroups of
$L_s$, $S \mapsto K_0(k[\eta^{-1}(S) \times \mathbb {Z}^n])$ is a
Frobenius module over Swan's Frobenius functor $S \mapsto G_0(\mathbb
{Z}S)$ \cite[\S 2]{Swan70}.  In view of important fact
\eqref{Eswan4}, to prove Theorem \ref{Tfiber0} it suffices
to show the following condition.
\begin{Cond} \label{Ccondition}
Some nonzero multiple of $x(S)$ is a sum of elements induced from
$K_0(kG)$ as $G$ varies over the finite subgroups of $\eta^{-1}(S)$.
\end{Cond}
Of course we only need to check Condition \ref{Ccondition} for one
group $S$ in each conjugacy class of cyclic subgroups of $L_s$.
Let $\sigma \colon L_s \to G$ denote the natural projection.  By
Remark \ref{Rconjugate} if $\sigma (S) = G$, then $S$ is conjugate to
$G_s$.  Hence we need only check Condition \ref{Ccondition} for $G_s$
and those cyclic subgroups $S$ such that $\sigma (S) \ne G$.

We start checking Condition \ref{Ccondition}
by considering the case where
$\sigma (S) \ne G$.  It clearly suffices to show that $\eta^{-1}(S)$
is a group for which Theorem \ref{Tfiber0} has already been verified
because it is lower in the lexicographic order.  To see this consider
the short exact sequence
\[
1 \longrightarrow \pi \longrightarrow \eta^{-1}(S) \longrightarrow
L(S) \longrightarrow 1
\]
where $L(S)$ denotes the inverse image of $S$ under the epimorphism
$L \twoheadrightarrow L_s$.  Observe that $L(S)$ is a
crystallographic group with
\begin{itemize}
\item $\dim L(S) = \dim (L)$ but
\item $\#(L(S)) < \#(L)$
\end{itemize}
since $\sigma (S) \ne G$.  Also notice that if $\hat {\pi}$ is any
subgroup of $\eta^{-1}(S)$ which contains $\pi$ with finite index,
then $\hat {\pi}$ satisfies conditions of Theorem \ref{Tfiber0}
because $\eta^{-1}(S) \subseteq \Gamma$ and we have assumed the same
property for $\Gamma$.  Since $\eta^{-1}(S)$ is lower in the
lexicographic order, we conclude from our inductive assumption that
$x(S)$ satisfies Condition \ref{Ccondition}.

It remains to show that $x(G_s)$ satisfies Condition \ref{Ccondition}.
Denote this element by $x_s$.  We will show by choosing the prime $s$
to be sufficiently large that $x_s$ also indeed satisfies Condition
\ref{Ccondition}.  Clearly we may assume that $m=0$.  Also we work
with the more convenient functor $\tilde{K}_0$ instead of $K_0$.
We can do this because of the following exact sequence
\[
0 \longrightarrow K_0(k)
\longrightarrow K_0(kH)
\longrightarrow \tilde{K}_0(kH)
\longrightarrow 0
\]
which is natural with respect to induction as $H$ varies over the
subgroups of $\eta^{-1}(S)$.

Now we assert that to show $x_s$ satisfies Condition
\ref{Ccondition}, it is enough to construct $s=s(x)$ so that $x_s$ is
in the image of the assembly map
\[
\mathfrak{a}_s \colon \mathbb {H}_0(\mathbb
{R}^n/L;\underline{\Wh}^k(q_s)) \longrightarrow
\tilde{K}_0(k\Gamma_s)
\]
(where, as in the proof of Theorem \ref{Tfiber}, $\Wh_1^k(\dummy) =
\Wh^k(\dummy)$, $\Wh_0^k(\dummy) = \tilde{K}_0(\dummy)$, and
$\Wh_j^k(\dummy) = K_j(\dummy)$ if $j<0$).  This assertion can be
seen as follows.  Since $L$ is a crystallographic group, $\mathbb
{R}^n$ has a triangulation on which $L$ acts simplicially and induces
a finite simplicial complex structure on $B = \mathbb
{R}^n/L$.  Let $B^i$ denote the $i$-skeleton of
$B$ and $E_s^i = q_s^{-1}(B^i)$.  Consider the
following diagram:
\begin{equation} \label{Equinn1}
\begin{CD}
E_s^0 @>\subseteq>>E_s @>\id>>E_s\\
@Vq_sVV @Vq_sVV @VVV\\
B^0 @>\subseteq>> B @>>> *.
\end{CD}
\end{equation}
The left commutative square in \eqref{Equinn1} induces a group
homomorphism
$\sigma \colon \mathbb {H}_0(B^0;\underline{\Wh}^k(q_s)) \to \mathbb
{H}_0(B;\underline{\Wh}^k(q_s))$ while the right square induces the
assembly map $\mathfrak{a}_s$ displayed above.  And $\sigma$ fits
into a (homology) exact sequence of abelian groups whose next group
is $\mathbb {H}_0(B,B_0;\underline{\Wh}^k(q_s))$ by \cite[Proposition
8.4]{Quinn82}.  This relative homology group can be analyzed by an
Atiyah-Hirzebruch-Quinn spectral sequence, constructed in
\cite[Proposition 8.7]{Quinn82}, which converges to $\mathbb
{H}_0(B,B_0; \underline{\Wh}^k(q_s))$ with $E^2_{j,-j} = H_j(B,B_0;
\Wh^k_{-j}(q_s))$ where $\Wh^k_{-j}(q_s)$ is the stratified system of
groups $\{\Wh^k_{-j}(\pi_1 q_s^{-1}(y))\mid y \in B\}$.  Now note that
\[
H_j(B,B_0;\Wh^k_{-j}(q_s)) =
\begin{cases}
0 &\text{if } j \notin[1,n]\\
0 \mod \mathcal{T} &\text{if } j>0.
\end{cases}
\]
The top equation is immediate for dimension reasons while the bottom
equation is a consequence of the fact that each $\pi_1(q_s^{-1}(y))$
contains $\pi$ with finite index.  Consequently $\mathbb {H}_0(B,B_0;
\underline{\Wh}^k(q_s)) \in \mathcal{T}$ and therefore also
$\coker(\sigma) \in \mathcal{T}$.  And consequently some nonzero
multiple of $x_s$ is in $\im(\mathfrak{a}_s \circ \sigma)$.  Next
consider the following second diagram:
\begin{equation} \label{Equinn2}
\begin{CD}
E_s^0 @>\id>> E_s^0 @>\subseteq>> E_s\\
@Vq_sVV @VVV @VVV\\
B^0 @>>> * @>>> *.
\end{CD}
\end{equation}
The left commutative square in \eqref{Equinn2} induces the assembly
map
\[
\mathfrak{a}_s^0 \colon \mathbb {H}_0(B^0; \underline{\Wh}^k(q_s))
\longrightarrow \bigoplus_{v \in B^0} \tilde{K}_0(k\pi_1
q_s^{-1}(v))
\]
while the right commutative square induces
\[
\tau \colon \bigoplus_{v\in B^0} \tilde{K}_0(k\pi_1q_s^{-1}(v))
\longrightarrow \tilde{K}_0(k\Gamma_s)
\]
which is the sum of the homomorphisms induced by the group inclusions
\[
\pi_1(q_s^{-1}(v)) \subseteq \pi_1(E_s) = \Gamma_s
\]
for $v \in B^0$.  Since the concatenation of the left and right
squares in \eqref{Equinn1} is the same square, namely
\[
\begin{CD}
E_s^0 @>\subseteq>>E_s \\
@Vq_sVV @VVV\\
B^0 @>>> *,
\end{CD}
\]
as the concatenation of the corresponding squares in \eqref{Equinn2},
we have by functoriality that $\mathfrak{a}_s \circ \sigma = \tau
\circ \bar{\mathfrak{a}}_s$ and consequently that
\begin{equation} \label{Equinn3}
\text{some nonzero multiple of $x_s$ is in $\im (\tau)$.}
\end{equation}
Recall that each group $\pi_1(q_s^{-1}(v))$,
$v \in B^0$, contains $\pi$
with finite index.  Hence the hypothesis of Theorem \ref{Tfiber0}
together with fact \eqref{Equinn3} verifies our assertion.

Now to show that $x_s \in \im(\mathfrak{a}_s)$,
it is easier to work with the Quinn assembly map at the
Whitehead group level; i.e.\ with the map
\[
\bar{\mathfrak{a}}_s \colon
\mathbb {H}_1 \bigl(\mathbb {R}^n/L; \underline {\Wh}^k
(\bar{q}_s) \bigr) \longrightarrow
\Wh^k(\Gamma \times \mathbb {Z}).
\]
where $\bar{q}_s \colon E_s \times S^1 \to \mathbb {R}^n/L$ is the
composite of the projection $E_s\times S^1 \to E_s$
with $q_s \colon E_s \to \mathbb {R}^n/L$.
In particular the following statement is true for $\bar
{\mathfrak{a}}_s$.
\begin{Ass} \label{Aassertion}
Given any element $y \in \Wh^k(\Gamma \times \mathbb {Z})$, there
exists a prime $\bar{s}(y) = s$ with $s \equiv 1 \mod | G
| $, such
that the transfer $y_s \in \Wh^k(\Gamma_s \times \mathbb {Z})$ of $y$
is in $\im \bar{\mathfrak {a}}_s$.
\end{Ass}
Before verifying Assertion \ref{Aassertion}, we use it to complete
the proof of Theorem \ref{Tfiber0}.  There is a natural ring
homomorphism $R[t,t^{-1}] \to \mu R$, where $R$ is any ring with 1
and $\mu R$ denotes its suspension, which induces the projection map
in the Bass-Heller-Swan formula $\underline{K}_*(R[t,t^{-1}]) \to
\underline{K}_{*+1}(R)$ on the spectrum level; cf.\ \cite{Gersten72},
\cite{Wagoner72}.  Consequently
there is the following commutative diagram
\[
\begin{CD}
\mathbb {H}_0 \bigl(\mathbb {R}^n/L; \underline{\Wh}^k
(q_s)\bigr)
@>\mathfrak{a}_s>>
\tilde{K}_0(k\Gamma_s) \\
@AAA @AAA \\
\mathbb {H}_0 \bigl(\mathbb {R}^n/L; \underline{\Wh}^k
(\bar{q}_s)\bigr)
@>\bar{\mathfrak{a}}_s>>
\Wh^k(\Gamma_s \times \mathbb {Z}).
\end{CD}
\]
Now set $s(x) = \bar{s}(y)$ where $y$ is the image of $x$ in
$\Wh^k(\Gamma \times \mathbb {Z})$ under the natural embedding
\[
\tilde{K}_0 (k\Gamma) \longrightarrow \Wh^k(\Gamma \times \mathbb
{Z}).
\]
Then by a diagram chase one sees directly that $x_{s(x)}$ is in the
image of $\mathfrak{a}_{s(x)}$.

It remains to verify Assertion \ref{Aassertion}.  For this purpose,
recall that Quinn showed that $y$ can be represented by a geometric
isomorphism (denoted $h$) of geometric $k$-modules on the space $E
\times S^1$ \cite[\S 3.2]{Quinn85x}.  The transfer of $h$ (denoted
$h_s$) to the finite sheeted cover
\[
E_s \times S^1 \longrightarrow E \times S^1
\]
clearly represents the transfer of
$y$ (denoted $y_s$) to $\Wh^k(\Gamma_s \times
\mathbb {Z})$.  But the radii $\epsilon_s$ of these geometric
isomorphisms $h_s$ measured via $\bar{q}_s$ in $\mathbb {R}^n/L$ go
to zero as $s \to \infty$; i.e.
\begin{equation} \label{Eepsilon0}
\lim_{s \to \infty} \epsilon_s = 0.
\end{equation}
In fact $h_s$ represents
an element $\bar {y}_s \in \Wh^k(\mathbb {R}^n/L, \bar{q}_s,
\epsilon_s)$ and $\bar{y}_s$ maps to $y_s$ under the natural
homomorphism
\[
\Wh^k(\mathbb{R}^n/L, \bar{q}_s, \epsilon_s) \longrightarrow
\Wh^k(\Gamma_s \times \mathbb {Z}).
\]
Since the strata in $\mathbb {R}^n/L$ of the stratified system
$\bar{q}_s$ are independent of $s$, Quinn's Stability Theorem
\cite[\S 4]{Quinn82} together with equation \eqref{Eepsilon0} and
\cite[\S 3]{Quinn85} yield that $y_s$ is contained in the image of
$\mathbb {H}_1(\mathbb {R}^n/L; \underline{\Wh}^k (\bar{q}_s))$,
i.e.\ $y_s \in \im (\bar{\mathfrak{a}}_s)$ when $s$ is sufficiently
large.  This verifies Assertion
\ref{Aassertion}, thus completing the proof of Theorem \ref{Tfiber0}
\end{proof}

We now use Theorem \ref{Tfiber0} to obtain the following result.
\begin{Cor} \label{Cfiber0}
Let $k$ be a field of prime characteristic and let
$\pi \lhd \Gamma$ be groups such that $\Gamma/\pi$
is an elementary amenable group.  Suppose that
$K_0(k[\tilde {\pi} \times \mathbb {Z}^m]) \otimes \mathbb {Q}$ is
generated by the images of $K_0(kG) \otimes \mathbb {Q}$ as $G$
varies over the finite subgroups of $\tilde {\pi}$ where $\tilde
{\pi}/ \pi$ is a finite subgroup of $\Gamma / \pi$,
for all $m \ge 0$.  Then for every $m \ge 0$,
\begin{enumerate} [\normalfont (i)] \label{Cfiber0i}
\item $K_i(k[\Gamma \times \mathbb {Z}^m]) \otimes \mathbb {Q}
= 0$ for all $i < 0$.
\item \label{Cfiber0ii}
$K_0(k[\Gamma \times \mathbb {Z}^m])
\otimes \mathbb {Q}$ is generated by the images
of $K_0(kG) \otimes \mathbb {Q}$ as $G$
varies over the finite subgroups of $\Gamma$.
\end{enumerate}
\end{Cor}
Note that Theorem \ref{TK_0}
immediately follows from using Corollary
\ref{Cfiber0}\eqref{Cfiber0ii}
in the case $\pi = 1$.  To see this, it will be
sufficient to show that $K_0(k[\tilde{\pi} \times \mathbb {Z}^m])
\otimes \mathbb {Q}$ is
generated  by the image of $K_0(k \tilde{\pi})
\otimes \mathbb {Q}$ whenever
$\tilde{\pi}$ is a finite subgroup of $\Gamma$, for all $m \ge 0$. 
But we know this by Lemma \ref{Lfinitegroup}.

\begin{proof}[Proof of Corollary \ref{Cfiber0}]
By Corollary \ref{Cnil4}, we need only prove
Corollary \ref{Cfiber0}\eqref{Cfiber0ii}.
We shall use the description of elementary amenable groups as
described in Lemma \ref{LEA}.
Let $\alpha$ be the least ordinal such that $\Gamma/\pi \in \mathcal
{X}_{\alpha}$.  Now $\alpha$ cannot be a limit ordinal, and the
result is clearly true if
$\Gamma/\pi \in \mathcal {X}_0$ because then
$\Gamma/\pi$ is finite.  Therefore we may assume that $\alpha$ is a
successor ordinal, and by transfinite induction that the result
is true for groups in $\mathcal {X}_{\alpha -1}$.
Thus $\Gamma$ has a normal subgroup $\pi_1$ containing
$\pi$ such that $\Gamma/\pi_1
\in \mathcal {A}$ and $\pi_1/\pi \in L\mathcal {X}_{\alpha - 1}$.
Now any $\mathcal {A}$-group maps onto a crystallographic group with
finite kernel; in other words there is a finite normal subgroup
$\pi_2/\pi_1$ of $\Gamma/\pi_1$ such that $\Gamma/\pi_2$ is
crystallographic.  Let
$\tilde {\pi}/\pi_2$ be any finite subgroup of $\Gamma/\pi_2$.
Then Lemma \ref{LEA} shows that $\tilde{\pi}/\pi
\in L\mathcal {X}_{\alpha -1}$.  Since $K_0$
commutes with direct limits,
Theorem \ref{Cfiber0}\eqref{Cfiber0ii} is true with $\tilde{\pi} =
\Gamma$ by our inductive hypothesis,
consequently $K_0(k[\tilde {\pi} \times \mathbb {Z}^m])$
is generated by the images of $K_0(kG) \otimes \mathbb {Q}$ as $G$
varies over the finite subgroups of $\tilde {\pi}$.
The result now follows from Theorem \ref{Tfiber0}.
\end{proof}

\section{The strong Bass conjecture and idempotents}
\label{SBass}

We begin this section by reviewing some results on lifting
idempotents modulo a nilpotent ideal, the relationship
between idempotents and projective modules, and the Hattori-Stallings
trace.  For more information, see \cite{Bass79} and \cite[\S III.7
and \S III.8]{Jacobson64}.  Some of the ideas used below originate
from \cite{Cliff80,FarkasMarciniak82}.

Let $R$ be a ring and let $d$ be a positive integer.
If $e \in \Mat_d(R)$ is an idempotent matrix, then $e$
defines by left multiplication an $R$-map $R^d \to R^d$ whose
image is a finitely generated projective
right $R$-module $P$.  Conversely given a finitely generated
projective $R$-module $P$, then for some positive
integer $d$ we may write $P \oplus Q \cong R^d$, and then the
projection from $R^d$ onto $P$ associated with this direct sum yields
an idempotent $e \in \Mat_d(R)$.  In this situation we shall say
that $P$ and $e$ correspond.  Of course $P$ does not determine $e$,
but if $e'$ is another choice for $e$,
then for $m$ sufficiently large, we may view $e, e' \in \Mat_m(R)$
with the property that $e' = ueu^{-1}$ for some $u \in \GL_m(R)$.

Suppose now that $N$ is a nil ideal of $R$ (so given $n \in N$, there
exists $s > 0$ such that $n^s = 0$).  If $\bar{e} \in
\Mat_d(R/N)$ is an idempotent, then
by the theory of \cite[\S III.8]{Jacobson64}, we may lift
$\bar{e}$ to an idempotent $e \in \Mat_d(R)$.
Since $N$ is contained in the Jacobson radical of $R$, by Nakayama's
lemma a map $P \to P'$ between finitely generated projective
$R$-modules is an isomorphism if and only
if the induced map $P/PN \to P'/P'N$ is an isomorphism.
Thus though $\bar{e}$ does not in general determine $e$, if $e'$ is
another choice for $e$, then the corresponding projective $R$-modules
$P,P'$ are isomorphic.

Now let $k$ be a commutative ring and let $\Gamma$ be a group.
If $\alpha \in \Mat_d(k\Gamma)$, then we may write $\alpha = \sum_{g
\in \Gamma} \alpha_g g$ where $\alpha_g \in \Mat_d(k)$.
Let $\tr \colon \Mat_d(k) \to k$ denote the usual trace map (i.e.\
the sum of the diagonal entries).  Then for $g \in \Gamma$
we define $\Tr_g \colon \Mat_d(k\Gamma) \to k$ by
\[
\Tr_g (\alpha) = \sum_{x \sim g} \tr (\alpha_x)
\]
where $x \sim g$ means that $x$ is conjugate to $g$ in $\Gamma$.
If $u \in \GL_d(k\Gamma)$, then
$\Tr_g(u \alpha u^{-1}) = \Tr_g(\alpha)$.
Suppose now $P$ is a finitely generated projective $k\Gamma$-module
and $P$ corresponds to an idempotent $e \in \Mat_d(k\Gamma)$.  Then we
set $r_P(g) = \Tr_g(e)$.  This is the Hattori-Stallings trace and
is well defined.  Moreover if $Q \cong P$ and $h \sim g$, then
$r_Q(h) = r_P(g)$.

Suppose $G$ is a group, $\theta \colon G \to \Gamma$ is a
homomorphism, and $Q$ is a finitely generated projective $kG$-module.
Then $Q \otimes_{kG} k\Gamma$ is a finitely generated projective
$k\Gamma$-module and we have the following induction formula:
\begin{equation} \label{Einduction1}
r_{Q \otimes_{kG} k\Gamma} (g) = \sum_{x \# g} r_Q(x),
\end{equation}
where $x \# g$ means that $\theta (x)$ is conjugate to $g$ in
$\Gamma$, and we choose only one $x$ from each $G$-conjugacy
class.  This means
in particular that if $g$ is not conjugate to any element of
$\theta (G)$, then $r_{Q \otimes_{kG} k\Gamma}(g) = 0$.

Now let $S$ be a complete commutative local ring with maximal ideal
$J$ such that $S/J \cong k$.
Suppose $P$ is a finitely generated projective $k\Gamma$-module.
Then $P$ corresponds to an idempotent $e \in \Mat_d(k\Gamma)$ for
some positive integer $d$.  Now
$k\Gamma \cong S\Gamma /J \Gamma $
and $J^n \Gamma /J^{n+1} \Gamma $ is a nil ideal
(in fact even the square of the ideal is 0) of
$S\Gamma / J^{n+1} \Gamma$ for all positive integers
$n$.  Set $e_1 = e$ and let
\[
\bar{\ } \colon
\Mat_d(S\Gamma /J^{n+1}\Gamma) \longrightarrow
\Mat_d(S\Gamma /J^n\Gamma)
\]
denote the natural epimorphism.
We shall also let $\bar{\ }$ denote the natural epimorphism
$S\Gamma/J^{n+1}\Gamma \twoheadrightarrow S\Gamma/J^n\Gamma$.
Thus starting with $n=1$, we may inductively lift $e_n \in
\Mat_d(S\Gamma /J^n \Gamma)$ to an idempotent $e_{n+1} \in
\Mat_d(S\Gamma / J^{n+1} \Gamma)$ such that
$\overline{e_{n+1}} = e_n$.  In particular
\[
\overline{\Tr_g(e_{n+1})} = \Tr_g(e_n)
\]
for all positive integers $n$ and for all $g \in \Gamma$.  This means
that for each $g$, the sequence $\bigl(\Tr_g(e_n) \bigr)$ yields a
well defined element of $S$, which we shall denote by $\hatTr_g(e)$.
Though the $e_n$ are not uniquely
determined by $e$, we do know that if $(e_n')$ is another such
sequence of idempotents, then
\[
\Tr_g(e_n') = \Tr_g(e_n)
\]
for all $n$, consequently $\hatTr_g(e)$ does not depend on the choice
of liftings for the $e_n$.  
Suppose $e' = ueu^{-1}$ where $u \in \GL_d(k\Gamma)$.  Then we
may lift the $u$ to units $u_n \in \GL_d(S\Gamma / J^n \Gamma)$ to
obtain a sequence of idempotents $u_n e_n u_n^{-1} \in \Mat_d(S\Gamma/
J^n \Gamma)$ which lift the idempotent $e'$.  This shows that
$\hatTr_g(e) = \hatTr_g(e')$.  Now
set $\hat{r}_P(g) = \hatTr_g(e)$.  Then the
above discussion shows that $\hat{r}_P(g)$ is well defined and if
$\bar{\ } \colon S\Gamma \to k\Gamma$ is the natural epimorphism,
then
\[
\overline{\hat{r}_P(g)} = r_P(g).
\]
Also if $U$ is a finitely generated projective $S\Gamma$-module such
that $U/UJ \cong P$, then $r_U(g) = \hat{r}_P(g)$.  For $\Gamma$
finite, such a $U$ will always exist
(this well-known fact can be seen from the above, because
the sequence of lifted idempotents $e_n$ will yield an idempotent
of $\Mat_d(S\Gamma)$, and we can let $U$ be the projective
$S\Gamma$-module corresponding to this idempotent; of course if
$\Gamma$ is infinite, this will not work because the supports of the
lifted idempotents $e_n$ may become arbitrarily large).
Finally if $\theta \colon G \to \Gamma$ is a group homomorphism and
$Q$ is a finitely generated projective $kG$-module, then from
\eqref{Einduction1} we get the corresponding  induction formula for
lifted traces:
\begin{equation}  \label{Einduction2}
\hat{r}_{Q \otimes_{kG} k\Gamma} (g) = \sum_{x \# g}
\hat {r}_Q(x).
\end{equation}

\begin{Rem} \label{RBrauer}
If in addition $p$ is a prime,
$k$ is a field of characteristic $p$, $S$ is a Noetherian
integral domain, and $\phi$ is the Brauer character of $P$,
then $\zeta (g) \hat{r}_P(g) =  \phi(g)$, where $\zeta (g)$ is the
order of the centralizer of $g$ in $\Gamma$.
\end{Rem}

Suppose $Q$ is a finitely generated
projective $k\Gamma$-module.  Then
$Q$ corresponds to an idempotent $f \in \Mat_c(k \Gamma)$ for some
positive integer $c$ and $e
\oplus f$ corresponds to $P\oplus Q$; here $e \oplus f$ means the
element in $\Mat_{d+c}(k \Gamma)$
\[
\begin{pmatrix}
e & 0 \\
0 & f
\end{pmatrix}.
\]
Thus if $(f_n)$ is a sequence of idempotents corresponding to $Q$,
then $(e_n \oplus f_n)$ is a sequence of idempotents corresponding to
$P\oplus Q$.  Since $\Tr_g(e_n \oplus f_n) = \Tr_g(e_n)+ \Tr_g(f_n)$
and $\hat{r}_P(g)$ is well defined, we
deduce that
\[
\hat{r}_{P \oplus Q}(g) = \hat{r}_P(g) + \hat{r}_Q(g).
\]
We can now prove (retaining the notation above)
\begin{Lem} \label{LBass}
Let $p$ be a prime, let $k$ be a field of characteristic $p$,
let $\Gamma$ be an elementary amenable
group, let $P$ be a finitely generated projective
$k\Gamma$-module, and let $g \in \Gamma$.  Suppose
all positive integers are nonzero divisors in $S$.
If $p | o(g)$ or $o(g) = \infty$, then
$\hat{r}_P(g) = 0$.  In particular, $r_P(g) = 0$.
\end{Lem}
\begin{proof}
By Theorem \ref{TK_0}, there is a positive integer $n$,
a finitely generated free $kG$-module $Q$,
finite subgroups $G_1, \dots G_m$ of $\Gamma$, and for each positive
integer $i \le m$ a finitely generated projective $kG_i$-module
$Q_i$ such that
\begin{align*}
P^n \oplus Q&\cong
\bigoplus_{i=1}^m Q_i \otimes_{kG_i} k\Gamma \oplus Q.\\
\intertext{Therefore}
n \hat{r}_P(g) = \sum_{i=1}^m \hat{r}_{Q_i \otimes_{kG_i} k\Gamma}
(g).
\end{align*}
Since $n$ is a nonzero divisor in $S$, it is sufficient to show that
$\hat{r}_{Q_i \otimes_{kG_i} k\Gamma} (g) = 0$ for all $i$.
Now $G_i$ is finite, so there exists a finitely generated projective
$SG_i$-module $\hat{Q}_i$ such
that $\hat{Q}_i/\hat{Q}_iJ \cong Q_i$.  Then applying
\eqref{Einduction2}, we have
\[
\hat{r}_{Q_i \otimes _{kG_i} k\Gamma} (g) = \sum_{x \# g}
r_{\hat{Q}_i} (x).
\]
In the case $o(g) = \infty$, this sum is empty and hence equal to
zero.  Also $r_{\hat{Q}_i} (x) = 0$ whenever $p | o(x)$
(this is really \cite[Theorem 59.7(ii)]{Dornhoff72} and Remark
\ref{RBrauer}).  Thus the result follows.
\end{proof}

\begin{Rem} \label{Rpadic}
Suppose $k$ is a finite field of characteristic $p$ and order $q$
where $q$ is a power of $p$.  Let $\zeta$ be a primitive
$(q-1)$th root of unity (in the algebraic closure of
$\mathbb {Q}_p$) and set $S = \mathbb {Z}_p[\zeta]$.
Then $S$ is a complete local ring with residue field $k$.
\end{Rem}

Using Lemma \ref{LBass}, we can now prove Theorem \ref{TBass}.
\begin{proof}[Proof of Theorem \ref{TBass}]
If $P$ corresponds to the idempotent $e$, then $e \in k'\Gamma$ for
some finitely generated subring $k'$ of $k$.  This means we may
assume that $k$ is finitely generated as a ring.
We have two cases to consider, depending on whether $k$ has
characteristic $p$ or characteristic $0$.
\bigskip

\noindent
Case 1.
\newline
$k$ has characteristic $p$.  Suppose by way of contradiction
that $g \in \Gamma$, $o(g)$ is not
invertible in $k$, yet $r_P(g) \ne 0$.  Then either $p | o(g)$
or $o(g) = \infty$.  Since $k$
is a finitely generated integral domain, it has a maximal ideal $M$
such that the image of $r_P(g)$ in $k/M$ is nonzero \cite[Ex.~24,
p.~71, Chapter 5]{AtiyahMacdonald69}.  Also $k/M$
will be a finite field
\cite[Ex.~6, p.~84, Chapter 7]{AtiyahMacdonald69}
of characteristic $p$.
Thus by Remark \ref{Rpadic}, there is a complete local ring
$S$ with residue
field $k/M$, and we can now obtain a contradiction by applying Lemma
\ref{LBass}.
\bigskip

\noindent
Case 2.
\newline
$k$ has characteristic 0.  Since $o(g)$ is not invertible in $k$, we
may choose a prime $p$ such that $pk \ne k$ and either $p | o(g)$
or $o(g) = \infty$.  Let $J$
be a maximal ideal in $k$ containing $p$, and let $S$ denote the
completion of $k$ with respect to $J$.  Then $S$ is a complete local
ring containing $k$
with residue field of characteristic $p$.  Also $k$ is a Noetherian
ring \cite[Corollary 7.7, p.~81]{AtiyahMacdonald69}, so by
\cite[Ex.~4, p.~114, Chapter 10]{AtiyahMacdonald69} it has the
property that all positive integers are nonzero divisors in $S$.
We now obtain a contradiction by applying Lemma \ref{LBass} and the
proof is complete.
\end{proof}

\begin{proof}[Proof of Theorem \ref{Tidemp}]
Suppose $e$ is a nontrivial idempotent in $k\Gamma$ and write $e =
\sum_{g \in \Gamma} e_gg$ where $e_g \in k$ for all $g \in
\Gamma$.  Since $e_g = 0$ for all but finitely many $g$, by replacing
$k$ with the subring generated by $\{e_g \mid e_g \ne 0\}$, we may
assume that $k$ is a finitely generated integral domain.  Choose $h
\in G\setminus 1$ such that $e_h \ne 0$, which is possible because
$e$ is nontrivial.  Then $k$ has a maximal ideal $M$ such that the
image of $e_h$ in $k/M$ is nonzero \cite[Ex.~24,
p.~71, Chapter 5]{AtiyahMacdonald69}.  Also $k/M$
will be a finite field
\cite[Ex.~6, p.~84, Chapter 7]{AtiyahMacdonald69}
of characteristic $p$.  Therefore we may assume that $k$ is a finite
field and hence by Remark \ref{Rpadic},
there is a complete local ring $S$ with residue
field $k$.  Now write $1 = e + f$ where $f$ is also an idempotent
in $k\Gamma$.  Let $P,Q$ be finitely generated
projective $k\Gamma$-modules which correspond to $e,f$
respectively, and
let $\theta \colon \Gamma \to 1$ denote the
natural epimorphism.  Then $P \otimes_{k\Gamma} k \oplus Q
\otimes_{k\Gamma} k \cong k$, so without loss of generality we may
assume that $P\otimes_{k\Gamma} k = 0$ and we see that $\hat{r} _{P
\otimes_{k\Gamma} k}(1) = 0$.
Since $\hat{r}_P(g) = 0$ for all $g \in \Gamma \setminus 1$ by
Lemma \ref{LBass}, we deduce from the induction formula
\eqref{Einduction2} that $\hat{r}_P(1) = 0$.
If $H$ is a finite subgroup of $\Gamma$, then $H$ is a
finite $p$-group, hence $kH$ is a local ring and so all
projective $kH$-modules are free.  Therefore $P^n$ is a stably free
$k\Gamma$-module for some positive integer $n$
by Theorem \ref{TK_0}, consequently $P^n \oplus (k\Gamma)^a \cong
(k\Gamma)^b$ for some integers $a,b$.  Since $\hat{r}_P(1) =0$, this
is only possible if $a=b$ and we now have $(k\Gamma)^a \cong
(k\Gamma)^a \oplus P^n$.

Since $e$ is nonzero by assumption, we see that $P^n$ is nonzero
(and also $a \ne 0$).
Therefore $(k\Gamma)^a$ is not a directly finite
$k\Gamma$-module \cite[\S
6.B]{Goodearl76}, hence by \cite[Lemma 6.9]{Goodearl76}
the ring $\Mat_a(k\Gamma)$ is not directly finite.  This means
that there exist $r,s \in \Mat_a(k\Gamma)$ such that $rs=1$ but $sr
\ne 1$.  By \cite[Theorem 3.2]{AraOmearaPerera}, all matrix rings
over the group algebra $kG$ of an amenable group $G$ are directly
finite.  Therefore $\Mat_a(k\Gamma)$ is directly finite, a
contradiction and the result is proven.
\end{proof}

\section{Whitehead groups of solvable linear groups} \label{Slinear}
In this section we shall prove Theorem \ref{Tlinear}.  The proof is
very similar to \cite[Theorem
1.1]{FarrellLinnell}, which showed that
$\Wh (G) = 0$ when $G$ is a virtually solvable linear group, so we
will not give full details.  For the rest of this section, we set
$k=k_2$.
We also need the following consequence of Theorem \ref{TFIC}(i).
\begin{Lem} \label{LFIC}
Let $H \lhd G$ be groups such that $G/H \in \mathcal {A}$, and let $p
\colon G \twoheadrightarrow G/H$ denote the natural epimorphism.
If $\Wh^k(p^{-1}(S) \times \mathbb {Z}^n) = 0$ for all virtually
cyclic subgroups $S$ of $G/H$ and for all non-negative integers $n$,
then $\Wh^k(G) = 0$.
\end{Lem}
\begin{proof}
This follows from Lemma \ref{Lliegroup} and Theorem \ref{TFIC}(i).
\end{proof}
\begin{proof}[Proof of Theorem \ref{Tlinear}]
Since Whitehead groups commute with direct limits, we may assume that
$G$ is finitely generated.  Also we may assume that $G > 1$.
We have two cases to consider.
\bigskip

\noindent
Case 1.  $G$ is virtually nilpotent.
\newline
Since $G$ is also finitely generated, we see that it is a virtually
poly-$\mathbb {Z}$ group.  Hence it is a discrete cocompact subgroup
of a virtually connected Lie group by a result of Auslander and
Johnson \cite[Theorem 1]{AuslanderJohnson76}.  And therefore $G$
satisfies the FIC in dimensions $\le 1$ because of Theorem
\ref{TFIC}.  Applying this to the situation where $\phi = \id \colon
\Gamma \to \Gamma$, we conclude that $\Wh^k(G) = 0$ since
\[
\Wh^k(C) = \tilde{K}_0(kC) = K_i(kC) = 0
\]
for all $i < 0$ when $C=1$ or $\mathbb {Z}$.  Since $\Wh^k(\dummy)$
commutes with direct limits and we have not used, so far, the
assumption that $G$ is a subgroup of $\GL_n(\mathbb {C})$, the
following additional remark is true.
\begin{Rem} \label{RWhnilpotent}
If $N$ is any torsion-free virtually nilpotent group, then $\Wh^k(N)
= 0$ and also $\tilde{K}_0(kN) = 0$.
\end{Rem}

To see the ``and also" assertion, notice that $\tilde{K}_0(kN)$ is a
subgroup of $\Wh^k(N\times \mathbb {Z})$
\cite[Chapter XII, \S 7]{Bass68} and that $N \times \mathbb
{Z}$ is also both torsion-free and virtually nilpotent.
\bigskip

\noindent
Case 2. $G$ is not virtually nilpotent.
\newline
By a theorem of Malcev
\cite[15.1.4]{Robinson96}, there exists $H \lhd G$
such that $H$ is nilpotent and $G/H \in \mathcal {A}$
(though $H$ will not be finitely generated in general).
\bigskip

Let $C/H$ be a virtually cyclic subgroup of $G/H$, let $n$ be a
non-negative integer, and set $D =  C \times \mathbb {Z}^n$.  By
Lemma \ref{LFIC}, it will suffice to show that $\Wh^k(D) = 0$.  Set
$E = H \times \mathbb {Z}^n$.  Then $E$ is a normal nilpotent
subgroup of $D$ and $D/E$ is virtually cyclic.  Since $D/E$ is
virtually cyclic, Proposition \ref{Pvircyclic} yields that
it has a finite normal subgroup $F/E$ such that $D/F$ is
either trivial, infinite cyclic or infinite dihedral.
For any finite subgroup $N/F$ of $D/F$ (e.g.\ $N = F$), the group $N$
is torsion-free virtually nilpotent and hence
\[
\Wh^k(N) = \tilde{K}_0(kN) = 0
\]
by Remark \ref{RWhnilpotent}.  In particular if $D/F = 1$, this shows
that $\Wh^k(D) = 1$.  In the other two cases, we now apply
Waldhausen's results \cite{Waldhausen78} and we shall
adopt his terminology.  We shall also require the following well known
result, which for example follows immediately from
\cite[Lemma 3.1]{FarrellLinnell}.
\begin{Lem} \label{Lregular}
Let $G$ be a finitely generated torsion-free virtually nilpotent
group.  Then $kG$ is regular Noetherian.
\end{Lem}
Using Lemma \ref{Lregular} and \cite[Theorem
19.1(iii)]{Waldhausen78}, we see that $kF$ is regular coherent.
If $D/F$ is infinite cyclic, then the result follows from
\cite[Corollary 17.2.3]{Waldhausen78}.  On the other hand if $D/F$ is
infinite dihedral, then the result follows from \cite[Corollary
17.1.3]{Waldhausen78}.
\end{proof}

\section{$k_2$ Group Rings} \label{Sappendix}

\begin{Thm} \label{TFIC}
The FIC in dimensions $\le 1$ is true for the group $\Gamma$
relative to the functor $\underline {K}(k_2\ )$ in the following
cases.
\begin{enumerate} [\normalfont(i)]
\item $\Gamma$ is a discrete cocompact subgroup of a
virtually connected Lie group.
\label{TFICi}
\item $\Gamma = H \wr \mathbb {Z}$ where $H$ is any finite
abelian group.
\label{TFICii}
\end{enumerate}
The precise statement is
that given any epimorphism $\phi \colon \pi \to \Gamma$, then
the assembly map
\[
\mathbb {H}_i \bigr(\mathcal {B}; \underline{K}(k_2
\Phi) \bigl) \longrightarrow K_i(k_2 \pi)
\]
induces an isomorphism for all $i \le 1$.  Here $\mathcal {B} =
\mathcal {E}/\Gamma$ where $\mathcal {E}$ is a universal
$\Gamma$-space with respect to the class of virtually cyclic
subgroups of $\Gamma$; $E$ is the universal cover of an
Eilenberg-Mac\,Lane space $K(\pi,1)$ and $\Phi \colon E \times_{\pi}
\mathcal {E} \to \mathcal{B}$ is induced by the projection $E \times
\mathcal {E} \to \mathcal {E}$.  Also $\pi$ acts on $\mathcal{E}$ via
$\phi$.
\end{Thm}

\begin{Cav}
The ``FIC in dimensions $\le 1$" defined in Theorem \ref{TFIC}
is somewhat weaker than the straightforward generalization of the FIC
which was formulated in \cite[\S 1.7]{FarrellJones93} and restated in
\cite[\S 7]{FarrellLinnell}.  In the notation of
\cite{FarrellLinnell}, the FIC allows an arbitrary free and properly
discontinuous $\Gamma$-space $Y$, while $Y = E/\ker \phi$ in Theorem
\ref{TFIC}.
\end{Cav}

Using Theorem \ref{TFIC}, we are able to prove Corollary \ref{CFIC},
which clearly implies Theorem \ref{Tfiber0} in the case $k = k_2$.
It now follows that we have a second proof of Theorem \ref{TK_0} in
the case $k = k_2$.  Also Theorem \ref{TFIC} enables us to prove 
Proposition \ref{Pwreath}.  On the other hand, we do not know how to
prove Theorem \ref{TFIC} for fields of characteristic different from
2; for example we do not know how to prove Theorem \ref{TFIC} in the
case $k = k_3$.  See \cite[Appendix]{FarrellJones93}
for details about universal
$\Gamma$-spaces relative to a class of subgroups of $\Gamma$.

\begin{Rem} \label{RFIC}
Formal homology properties show that Theorem \ref{TFIC} is equivalent
to the statement that the assembly map
\[
\mathbb{H}_i \bigl(\mathcal {B}; \underline{\Wh}^k(\Phi) \bigr)
\longrightarrow \Wh_i^k(\pi)
\]
induces an isomorphism for all $i \le 1$.  Here $\underline {\Wh}^k(
\pi)$ is the cofiber of the standard map of spectra
$\mathbb {H}(B\pi; \underline{K}(k)) \to \underline {K} (k\pi)$
defined by Loday \cite{Loday76} (which is the Quinn assembly map
associated to the fibration $\id \colon B\pi \to B\pi$; cf.\
\cite{Quinn85}).  Here $k = k_2$.
\end{Rem}

\begin{proof}[Proof of Proposition \ref{Pwreath}]
Since $\tilde{K}_0$ commutes with direct limits, we may assume that
$G$ is a finite abelian 2-group.
In view of Theorem \ref{TFIC}(ii), it will
now be sufficient to show that
$\tilde{K}_0(k_2 H) = 0$ and $K_i (k_2 H) = 0$ for every
virtually cyclic subgroup $H$ of $G \wr \mathbb {Z}$ and for every
$i < 0$.  But the only virtually cyclic subgroups of
$G \wr \mathbb {Z}$ are $\mathbb {Z}$ and finite
abelian 2-groups.  Therefore by Lemma
\ref{Lfinitegroup}, we are reduced to showing that $\tilde{K}_0(k_2H)
= 0$ when $H$ is a finite 2-group.  But this is obvious because now
$k_2H$ is a local ring and so all projective $k_2H$-modules are free.
\end{proof}

\begin{Cor} \label{CFIC}
Let $\pi \lhd \Gamma$ be groups with $\Gamma/\pi \in \mathcal {A}$.
Suppose for every $n \ge 0$ and every subgroup
$\tilde {\pi}$ of $\Gamma$ containing $\pi$ as a subgroup of finite
index, $K_0(k_2[\tilde{\pi} \times \mathbb {Z}^n])
\otimes \mathbb {Q}$ is generated by the images of $K_0(k_2G)$
as $G$ varies over the finite subgroups of $\tilde {\pi}$.
Then for every $n \ge 0$,
\begin{enumerate} [\normalfont (i)]
\item $K_i(k_2[\Gamma \times \mathbb {Z}^n])
\otimes \mathbb {Q} = 0$ for all $i < 0$.
\item $K_0(k_2[\Gamma \times \mathbb {Z}^n]) \otimes \mathbb {Q}$ is
generated by the images of $K_0(k_2G) \otimes \mathbb {Q}$ as $G$
varies over the finite subgroups of $\Gamma$.
\end{enumerate}
\end{Cor}
\begin{proof}
Let $L = \Gamma/\pi$, let $\phi \colon \Gamma \to L$ denote the
natural epimorphism, and set $k = k_2$.
Lemma \ref{Lliegroup} tells us that
$L$ is a discrete cocompact subgroup of a virtually connected
Lie group, hence Theorem \ref{TFIC}(i) is applicable.  We deduce that
it suffices to prove Corollary \ref{CFIC} under the extra assumption
that $L$ is either the infinite cyclic group $\mathbb {Z}$ or the
infinite dihedral group $C_2 * C_2$.  To see this reduction use
Corollary \ref{Cnil4} and the following facts:
\begin{itemize}
\item
$\mathbb {H}_i \bigl(\mathcal {B}; \underline{K}(k \Phi
)\bigr) = \varinjlim_{B} \mathbb {H}_i
\bigl (B; \underline {K} (k \Phi) \bigr)$ where $B$
varies over the finite subcomplexes of $\mathcal {B}$.

\item Every infinite virtually cyclic subgroup of $L$ maps onto either
$\mathbb {Z}$ or $C_2 *C_2$ with finite kernel.

\item There is a spectral sequence due to Atiyah-Hirzebruch-Quinn with
$E^2_{ij} = H_i \bigl(B; K_j(k \Phi)\bigr)$ converging to
$\mathbb {H}_{i+j} \bigl (B; \underline {K}(k \Phi)
\bigr)$ where $K_j(k\Phi)$ is the stratified system of groups
$\{K_j(k\phi^{-1}(G_x)) \mid x \in B\}$ and
each $G_x$ is a virtually cyclic group.
\end{itemize}
So we have two cases to consider.
\bigskip

\noindent
Case 1. $L=\mathbb {Z}$.
\newline
Then $\Gamma = \pi \rtimes \mathbb {Z}$ and the result is
an immediate consequence of Corollary \ref{Cnil4}.
\bigskip

\noindent
Case 2. $L = C_2 * C_2$.
\newline
This case would follow from an argument similar to that given for
Case 1 using \cite{Waldhausen78} in place of \cite{FarrellHsiang70}
provided we could prove an analogue of Lemma \ref{Lkey} for the
``Nil-groups" which occur in Waldhausen's formula.  But we suspect
that these more esoteric groups are ``dominated" by those already
analyzed.  (Some evidence for this suspicion is given by
\cite[Theorem 2.6]{FarrellJones95}.)
Hence we have opted instead for the
following argument which reduces Case 2 also to Lemma \ref{Lkey}.  We
proceed via a Frobenius induction argument which is similar to the
ones in Theorems \ref{Tfiber} and \ref{Tfiber0}
(but more explicit), and our notation
will be consistent with that used there.  Thus we have
\[
q \colon E \longrightarrow \mathbb {R}/ C_2 * C_2 = [0,1]
\]
with the action of $C_2 * C_2$ on $\mathbb {R}$ generated by the two
reflections: $a$ in 0 and $b$ in 1.  Note that $s$ can be any odd
prime in this case, since $\# (L) = 2$.  The endomorphism $\psi$
of $C_2 * C_2$ is explicitly described by the formula
\[
\psi (a) = a \quad \text{and} \quad \psi(b) = (ba)^{s-1} b.
\]
Likewise $f \colon \mathbb {R} \to \mathbb {R}$ is given by
\[
f(x) = sx \quad \text{for all } x \in \mathbb {R}.
\]
The quotient group $L_s$ of $L$ is now the finite dihedral group
$D_{2s}$ of order $2s$; it is gotten from $C_2 * C_2$ by adding the
extra relation $(ab)^s = 1$.  Note that $\psi (L)$ is the subgroup of
index $s$ in $C_2 * C_2$ generated by $a$ and $(ba)^{s-1}b$, and is
of course isomorphic to $C_2 * C_2$ via $\psi$.  Also $\Gamma_s =
\phi^{-1}(\psi (L))$ is a subgroup of index $s$ in $\Gamma$.
The smooth 5-manifold
$E_s$ is the space $\mathbb {R} \times_{\Gamma_s} \tilde
{M}$ and the map $q_s \colon E_s \to [0,1]$ is the composition of
\begin{align*}
E_s &\longrightarrow \mathbb {R}/\psi(K) = [0,s] \\
\intertext{with the homeomorphism}
\hat{f}^{-1} &\colon [0,s] \longrightarrow [0,1]
\end{align*}
given by multiplication by $1/s$.  It is a stratified system of
fibrations with $\{0,1\}$ and $(0,1)$ the strata for $[0,1]$.  We
will apply Frobenius induction with respect to the finite quotient
group $L_s$ of $\Gamma$.  Since we only need to obtain results modulo
the class $\mathcal {T}$, it suffices to consider only the maximal
cyclic subgroups of $L_s$ and we need only one representative from
each conjugacy class.  These can be chosen to be the cyclic subgroups
$\langle a \rangle$ and $\langle ab \rangle$.  Note that $|
\langle a \rangle | = 2$ and $| \langle ab \rangle
| = s$.  Let $\eta \colon
\Gamma \to L_s$ denote the composite epimorphism
\[
\Gamma \overset{\phi}{\twoheadrightarrow} L
\twoheadrightarrow L_s.
\]  
We will be concerned with the subgroups $\eta^{-1} (\langle a
\rangle)$ and $\eta^{-1}(\langle ab \rangle )$.  Note that $\eta^{-1}
(\langle a \rangle ) = \Gamma_s$ and $\phi \eta^{-1}
(\langle ab \rangle) \cong \mathbb {Z}$.
Denote $\eta^{-1}(\langle ab \rangle )$ by $\Gamma_{(s)}$.
Hence the already verified Case 1 yields that $K_0(k[\Gamma_{(s)}
\times \mathbb {Z}^n]) \otimes \mathbb {Q}$ is generated by the
images of $K_0(kG) \otimes \mathbb {Q}$ as $G$ varies over
the finite subgroups of $\Gamma_{(s)}$.
Thus by Corollary \ref{Cnil4} and the important fact \eqref{Eswan4},
to complete the proof in Case 2, it suffices
to associate an odd prime $s = s(x)$ to each element $x \in K_0(k
[\Gamma \times \mathbb {Z}^n])$ where $n \ge 0$, such that the
transfer $x_s$ of $x$ to $K_0(k[\Gamma_s \times \mathbb {Z}^n])$
satisfies the following:
\begin{quote}
There exists a positive integer $m$ such that $m x_s$ is the sum
of elements induced from $K_0(kG)$ as $G$ varies over the finite
subgroups of $\Gamma_s$.
\end{quote}
We may assume that $n = 0$ because the $\mathbb {Z}^n$ factor is a
``dummy variable", and we work with the more convenient functor
$\tilde {K}_0$ instead of $K_0$; this is clearly equivalent
(cf.\ Remark \ref{RFIC}).

Notice that it is enough to construct $s = s(x)$ so that $x_s$ is in
the image of the assembly map
\[
\mathfrak {a}_s \colon \mathbb {H}_0 \bigl([0,1]; \Wh^k
(q_s) \bigr) \longrightarrow \tilde {K}_0 (k
\Gamma_s).
\]
To see this note that the Atiyah-Hirzebruch-Quinn spectral sequence
together with our assumption of Corollary \ref{CFIC} yields that any
element in $\im (\mathfrak {a}_s)$ has a nonzero multiple which is
the sum of elements induced from $\tilde {K}_0(kG)$ as $G$ varies
over the finite subgroups of $\Gamma$.  (See the second paragraph
before Assertion 4.3 for more details.)

On the other hand, it was shown in the proof of Theorem \ref{Tfiber}
that an analogous statement with $\tilde {K}_0$ replaced by $\Wh^k$ is
true.  To be explicit given any $y \in \Wh^k(\Gamma \times \mathbb
{Z})$, there is an odd prime $\bar {s}(y) = s$ such that the transfer
$y_s \in \Wh^k (\Gamma_s \times \mathbb {Z})$ of $y$ is in the image
of the assembly map
\[
\bar {\mathfrak {a}}_s \colon \mathbb {H}_1([0,1]; \underline {\Wh}^k
(\bar{q}_s)) \longrightarrow
\Wh^k(\Gamma_s \times \mathbb {Z}).
\]
where $\bar{q}_s \colon E_s \times S^1 \to [0,1]$ is the composite of
the projection $E_s \times S^1 \to E_s$ with $q_s \colon E_s \to
[0,1]$.
Since the projection map from $\Wh$ to $\tilde {K}_0$ in the
Bass-Heller-Swan formula can be naturally defined on the spectrum
level, there is the following commutative diagram
\[
\begin{CD}
\mathbb {H}_0 \left([0,1]; \underline{\Wh}^k(q_s)
\right)
@>\mathfrak{a}_s>>
\tilde{K}_0 (k\Gamma_s) \\
@AAA @AAA \\
\mathbb {H}_1 \left([0,1]; \underline {\Wh}^k (
\bar{q}_s) \right)
@>\bar {\mathfrak {a}}_s>>
\Wh^k(\Gamma_s \times \mathbb {Z}).
\end{CD}
\]
Now set $s(x) = \bar{s}(y)$ where $y$ is the image of $x$ in
$\Wh^k(\Gamma \times \mathbb {Z})$ under the natural embedding
\[
\tilde {K}_0 (k\Gamma) \longrightarrow \Wh^k(\Gamma \times \mathbb
{Z}).
\]
Then by a diagram chase, one sees directly that $x_{s(x)}$ is in the
image of $\mathfrak {a}_{s(x)}$.
\end{proof}

\begin{proof}[Proof of Theorem \ref{TFIC}]
We start by making some general constructions of spectra arising from
a unital ring homomorphism by applying the $K$-theory functor.
\begin{Not*}
We let $f \colon R \to S$ denote a unital ring homomorphism between
rings $R$ and $S$ both of which contain 1, and we let $f_* \colon
\underline{K}(R) \to \underline{K}(S)$ denote the map of $K$-theory
spectra induced by $f$ with $\underline{K}(f)$ denoting the cofiber
spectrum of $f_*$.  If $G$ is a group, then $f^G \colon RG \to SG$
denotes the ring homomorphism induced by $f$.  Note that
\[
f_*^G \colon \underline {K}(RG) \longrightarrow \underline {K}(SG)
\]
denotes the corresponding map of spectra whose cofiber is denoted by
$\underline{K}(f^G)$.  Composition with the fundamental group
determines three spectra valued functors $F_R$, $F_S$, $F_f$ from
topological spaces $X$ to spectra explicitly given by the formulae
\[
F_R(X) = \underline {K}(R\pi_1(X)), \quad
F_S(X) = \underline {K}(S\pi_1(X)), \quad
F_f(X) = \underline {K}(f^{\pi_1(X)})
\]
and a natural transformation
\[
f^{\pi_1(X)}_* \colon F_R(X) \longrightarrow F_S(X).
\]
\end{Not*}
Let $p \colon E \to B$ denote a simplicially stratified fibration
(cf.\ \cite[p.~254]{FarrellJones93}).
Then $p$ determines three new spectra by applying Quinn's
\emph{Homology with stratified coefficients} construction to $p$
(cf.\ \cite[Appendix]{Quinn82}, \cite[\S 1.2
and \S 1.4]{FarrellJones93}), namely
\[
\mathbb {H}(B;F_R(p)), \quad
\mathbb {H}(B;F_S(p)), \quad
\mathbb {H}(B;F_f(p)), \quad
\]
and the natural transformation determines a (coefficient) map of
spectra
\[
f_B \colon \mathbb {H}(B;F_R(p)) \longrightarrow \mathbb {H}(B;
F_S(p))
\]
whose cofiber spectrum is $\mathbb {H}(B; F_f(p))$.  This fact
follows from the explicit construction of $\mathbb {H}(\mathcal{B};
\mathcal {S}(p))$ given in \cite[Definition 8.1]{Quinn82} (where $X
\mapsto \mathcal{S}(X)$ is an arbitrary covariant functor from spaces
to spectra) using that the cofiber spectra are defined by a mapping
cone construction; cf.\ \cite[p.~154]{Adams74}.  And Quinn's
construction \cite[Definition 8.1]{Quinn82} is applicable because of
\cite[Proposition 8.4]{Quinn82}.

For any functor $F$ from spaces to spectra, let
\[
\mathfrak {a}_F \colon \mathbb {H}(B; F(p)) \longrightarrow F(E)
\]
denote the Quinn assembly map (of spectra); cf.\
\cite[Appendix]{Quinn82}.  Since assembly respects natural
transformations of functors, we have the following
two step (homotopy) commutative ladder of spectra:
\begin{equation} \label{ETFIC1}
\begin{CD}
\underline{K}(R\pi_1(E))
@>f_*^{\pi_1(E)}>>
\underline{K}(S\pi_1(E))
@>>>
\underline{K}(f^{\pi_1(E)}) \\
@A\mathfrak{a}_{F_R}AA
@A\mathfrak{a}_{F_S}AA
@A\mathfrak{a}_{F_f}AA \\
\mathbb{H} (B; F_R(p))
@>f_B>>
\mathbb{H}(B; F_S(p))
@>>>
\mathbb {H}(B;F_f(p)).
\end{CD}
\end{equation}
Ladder \eqref{ETFIC1} yields two long exact sequences in homotopy
and a commutative ladder connecting them.

Recall there is a fiber square of unital ring epimorphisms
\begin{equation}  \label{ETFIC2}
\begin{CD}
\mathbb {Z}C_2 @>f_1>> \mathbb {Z} \\
@Vg_1VV @VVg_2V  \\
\mathbb {Z} @>>f_2> k_2
\end{CD}
\end{equation}
where $f_1$ and $g_1$ send the generator of $C_2$ to 1 and $-1$
respectively.  This square induces
the fiber square
\begin{equation} \label{ETFIC3}
\begin{CD}
(\mathbb{Z} C_2)\pi @>f_1^{\pi}>> \mathbb {Z} \pi \\
@Vg_1^{\pi}VV  @VVg_2^{\pi}V  \\
\mathbb {Z} \pi @>>f_2^{\pi}> k_2 \pi.
\end{CD}
\end{equation}
Let the simplicially stratified fibration $p$ be the projection map
\begin{equation} \label{Ep}
p \colon
\mathcal {E} \times_{\Gamma} \hat{B} \pi \longrightarrow \mathcal {B}
= \mathcal {E}/\Gamma
\end{equation}
where $\hat{B} \pi \to B\pi$ is the covering space of the
Eilenberg-Mac\,Lane space $B\pi$ for the group $\pi$ corresponding to
the normal subgroup $\ker \phi$ of $\pi$.

Now consider the commutative ladder of homology exact sequences
determined by ladder \eqref{ETFIC1} in the special case where $f
\colon R \to S$ is $f_1 \colon \mathbb {Z}C_2 \to \mathbb {Z}$ (and
$p \colon E \to B$ is as specified in \eqref{Ep}).  Note
that $\mathfrak {a}_{F_{\mathbb {Z}}}$ and $\mathfrak {a}_{F_{\mathbb
{Z}C_2}}$ are the assembly maps in the FIC for $\Gamma$ relative to
the functor $\underline {K}(\mathbb {Z} \ )$ and the epimorphisms
$\phi \colon \pi \to \Gamma$ and $\hat {\phi} \colon C_2 \times \pi
\to \Gamma$, respectively, where $\hat{\phi}$ is the composition of
$\phi$ with the canonical epimorphism
$C_2 \times \pi \twoheadrightarrow \pi$.  But these assembly maps
induce isomorphisms in homotopy in all dimensions $\le 1$, by
\cite[Theorem 2.1]{FarrellJones93} for Theorem
\ref{TFIC}\eqref{TFICi} and by \cite[Corollary 4.2]{FarrellLinnell}
for Theorem \ref{TFIC}\eqref{TFICii}.
(Here are some more details concerning the first assertion in the
last sentence.  Remark \ref{RFIC}, replacing $k$ by $\mathbb{Z}$,
reduces it to the same assertion with the spectrum $\underline{K}$
replaced by the spectrum $\underline{\Wh}$.  But there is a natural
map of spectra
\[
\psi_2\colon \underline{\mathcal {P}}_*(\dummy) \longrightarrow
\underline{\Wh}_{*-2}(\dummy)
\]
where $\underline{\mathcal{P}}(\dummy)$ is the topological
pseudo-isotopy spectrum.  And $\psi_2$ induces an isomorphism of
$\pi_i$ for all $i \le -1$.  In particular $\Wh(\pi_1 X) =
\pi_{-1}\underline{P}(X)$.  See
\cite[p. 775, Proof of Corollary 0.7]{FarrellJones91}
for a detailed account of an analogous assertion.)
Consequently so does $\mathfrak
{a}_{F_{f_1}}$.  We will soon show in Lemma \ref{LTFIC}
that this is also true for $\mathfrak{a}_{F_{f_2}}$.
To do this, first note that the fiber
square \eqref{ETFIC2} determines a natural transformation
\[
g^X \colon F_{f_1}(X) \longrightarrow F_{f_2}(X)
\]
making the following ladder of spectra maps commutative
\[
\begin{CD}
F_{\mathbb {Z}C_2}(X)  @>{f_1^{\pi_1(X)}}_*>>
F_{\mathbb {Z}}(X) @>>> F_{f_1}(X)  \\
@V{g_1^{\pi_1(X)}}_*VV @V{g_2^{\pi_1(X)}}_*VV
@Vg^XVV \\
F_{\mathbb {Z}}(X) @>>{f_2^{\pi_1(X)}}_*> F_{k_2}(X)
@>>> F_{f_2}(X).
\end{CD}
\]
The map of ``coefficients" induced by this natural
transformation is denoted by
\[
\hat{g} \colon \mathbb {H}(B; F_{f_1}(p)) \longrightarrow \mathbb
{H}(B; F_{f_2}(p))
\]
and fits into the following (homotopy) commutative diagram of spectra
maps:
\[
\begin{CD}
\underline{K} (f_1^{\pi_1(E)})  @>g^E>>
\underline{K}(f_2^{\pi_1(E)})  \\
@A\mathfrak {a}_{F_{f_1}}AA
@AA\mathfrak {a}_{F_{f_2}}A  \\
\mathbb {H}(B; F_{f_1}(p))
@>>\hat {g}>
\mathbb {H}(B; F_{f_2}(p)).
\end{CD}
\]
\begin{Lem} \label{LTFIC}
In this square both $g^E$ and $\hat {g}$ induce isomorphisms in
homology in all dimensions $\le 1$.
\end{Lem}
Before proving Lemma \ref{LTFIC}, note that it directly implies the
fact about $\mathfrak{a}_{F_{f_2}}$ asserted above; i.e.\ we have the
following consequence
\begin{Cor} \label{CTFIC}
The map $\mathfrak {a}_{F_{f_2}}$ (as well as $\mathfrak
{a}_{F_{f_1}}$) of spectra induces an isomorphism in homology in all
dimensions $\le 1$.
\end{Cor}
\begin{proof}[Proof of Lemma \ref{LTFIC}]
We start by establishing the general fact that $g^X$ induces an
isomorphism in homology in all dimensions $\le 1$ for every path
connected space $X$, and for $E$ in particular.  But this is merely a
reformulation of Bass's excision property which states that the
relative $K$-groups $K_i(R,I)$ for a two-sided ideal $I$ in a ring
$R$ depend, for all $i \le 0$, only on the ring structure of $I$ (not
on $R$).  Then note that
\begin{align*}
K_i(f_1^{\pi}) &= K_{i-1}(\mathbb {Z}[C_2 \times \pi], \ker
(f_1^{\pi})) \quad \text{and} \\
K_i(f_2^{\pi}) &= K_{i-1}(\mathbb {Z}\pi,
\ker(f_2^{\pi}))
\end{align*}
where $\pi = \pi_1(X)$.  But
$\ker (f_1^{\pi}) \cong \ker (f_2^{\pi})$ as rings, because
\eqref{ETFIC3} is a fiber square of epimorphisms.

To establish the statement about $\hat {g}$, we use a variant of the
comparison theorem for the map induced by $\hat {g}$ between the
Atiyah-Hirzebruch-Quinn spectral sequences for $\mathbb {H}(B;
F_{f_1}(p))$ and $\mathbb {H}(B; F_{f_2}(p))$, respectively.  By the
general fact about $g^X$ demonstrated above, we get that this map is
an isomorphism between the $E^2_{ij}$-terms for all $j \le 1$.  From
this one shows (by induction on $n$) that it induces an isomorphism
between $E^n_{ij}$-terms for all $i+j \le 1$ and an epimorphism
between $E^n_{i,2-i}$-terms for all $i \ge 1$.  Consequently
$\hat{g}$ induces an isomorphism in homology in all dimensions $\le
1$.
\end{proof}

Now consider the commutative ladder of homology exact sequences
determined by ladder \eqref{ETFIC1} in the new special case where $f
\colon R \to S$ is $f_2 \colon \mathbb {Z} \to k_2$ (and $p \colon E
\to B$ is as in \eqref{Ep}).  Note that
$\mathfrak {a}_{F_{k_2}}$ is the assembly map mentioned in Theorem
\ref{TFIC}.  Since $\mathfrak {a}_{F_{\mathbb {Z}}}$ and $\mathfrak
{a}_{F_{f_2}}$ both induce isomorphisms in homology in all dimensions
$ \le 1$ because of \cite[Theorem 2.1]{FarrellJones93},
\cite[Corollary 4.3]{FarrellLinnell} and Corollary \ref{CTFIC}, 
we conclude from the 5-lemma that
$\mathfrak {a}_{F_{k_2}}$ induces an isomorphism in homology in all
dimensions $\le 0$ and an epimorphism in dimension 1.

To complete the proof of Theorem \ref{TFIC}, it remains to show that
$\mathfrak{a}_{F_{k_2}}$ induces a monomorphism in dimension 1.  To
do this we make use of Milnor's Mayer-Vietoris type sequence in
$K$-theory; cf.\ \cite[Theorem 3.3]{Milnor71}.
We see from it that given an element $x \in K_1(\mathbb {Z}\pi)$
which maps to 0 in $K_1(k_2\pi)$ via the map induced by $f_2^{\pi}$;
i.e.\ $K_1(f_2^{\pi})(x) = 0$, then there exists an element $\hat {x}
\in K_1(\mathbb {Z}[C_2 \times \pi])$ satisfying

\begin{equation} \label{ETFIC4}
\begin{aligned}[c]
&1.\\
&2.
\end{aligned}
\quad
\begin{aligned}[c]
K_1(g_1^{\pi})(\hat{x}) &= x \\
K_1(f_1^{\pi})(\hat{x}) &= 0.
\end{aligned}
\end{equation}
\begin{Not*}
Let $A_{C_2}$, $A_{\mathbb {Z}}$ and $A_k$ denote the group
homomorphisms induced by the assembly maps $\mathfrak {a}_{F_{\mathbb
{Z}C_2}}$, $\mathfrak {a}_{F_{\mathbb {Z}}}$ and $\mathfrak
{a}_{F_{k_2}}$, respectively, in 1-dimensional homology.  Also let
$\bar{f}_1$, $\bar{f}_2$, $\bar{g}_1$ and $\bar{g}_2$ denote the
homomorphisms in 1-dimensional homology induced by the spectra maps
$(f_1)_B$, $(f_2)_B$, $(g_1)_B$ and $(g_2)_B$, respectively.
\end{Not*}

Let $z \in \ker (A_k)$.  Once we show that $z = 0$, we will have
completed the proof of Theorem \ref{TFIC}.  By the second commutative
ladder of homology exact sequences described above, $z$ pulls back to
an element $y \in \mathbb {H}_1(B; F_{\mathbb {Z}}(p))$; i.e.\
$\bar{f}_2(y) = z$.  Now let $x = A_{\mathbb {Z}}(y)$ and note that
$K_1(f_2^{\pi})(x) = 0$.  Hence there is an element $\hat {x} \in
K_1(\mathbb {Z}[C_2 \times \pi ])$ solving equations 1 and 2 of
\eqref{ETFIC4}.  Since $A_{C_2}$ is an isomorphism, there exists
\[
\hat{y} \in \mathbb {H}_1(B; F_{\mathbb {Z}C_2}(p))
\]
such that $A_{C_2}(\hat{y}) = \hat{x}$.  Note that equations 1 and 2
of \eqref{ETFIC4} yield
\[
\begin{aligned}[c]
&1.\\
&2.
\end{aligned}
\quad
\begin{aligned}[c]
\bar{g}_1(\hat{y}) &= y \\
\bar{f}_1 (\hat {y}) &= 0
\end{aligned}
\]
since $A_{\mathbb {Z}}$ is an isomorphism and the following two
diagrams commute:
\[\hspace*{-.8ex}
\begin{CD}
\mathbb {H}_1(B; F_{\mathbb {Z}C_2}(p))
@>\bar{g}_1>>
\mathbb {H}_1(B; F_{\mathbb {Z}}(p)) \\
@VA_{C_2}VV  @VVA_{\mathbb {Z}}V  \\
K_1(\mathbb {Z}[C_2 \times \pi ])
@>>K_1(g_1^{\pi})>
K_1(\mathbb {Z}\pi)
\end{CD}
\quad\text{and}\quad
\begin{CD}
\mathbb {H}_1(B; F_{\mathbb {Z}C_2}(p))
@>\bar{f}_1>>
\mathbb {H}_1(B; F_{\mathbb {Z}}(p)) \\
@VA_{C_2}VV  @VVA_{\mathbb {Z}}V  \\
K_1(\mathbb {Z}[C_2 \times \pi ])
@>>K_1(f_1^{\pi})>
K_1(\mathbb {Z}\pi).
\end{CD}
\]
But $\bar{g}_2 \circ \bar{f}_1 = \bar{f}_2 \circ \bar{g}_1$ because
of square \eqref{ETFIC2}, hence
\[
0 = \bar{g}_2 (\bar{f}_1 (\hat{y}))
= \bar{f}_2(\bar{g}_1(\hat{y})) = \bar{f}_2(y) = z.
\]
\end{proof}

\bibliographystyle{plain}

\begin{thebibliography}{10}

\bibitem{Adams74}
J.~F. Adams.
\newblock {\em Stable homotopy and generalised homology}.
\newblock University of Chicago Press, Chicago, Ill., 1974.
\newblock Chicago Lectures in Mathematics.

\bibitem{AraOmearaPerera}
Pere Ara, Kevin~C. O'Meara, and Francesc Perera.
\newblock Stable finiteness of group rings in arbitrary characteristic.
\newblock {\em Adv. Math.}, 170(2):224--238, 2002.

\bibitem{AtiyahMacdonald69}
M.~F. Atiyah and I.~G. Macdonald.
\newblock {\em Introduction to commutative algebra}.
\newblock Addison-Wesley Publishing Co., Reading, Mass.-London-Don Mills, Ont.,
  1969.

\bibitem{AuslanderJohnson76}
L.~Auslander and F.~E.~A. Johnson.
\newblock On a conjecture of {C}. {T}. {C}. {W}all.
\newblock {\em J. London Math. Soc. (2)}, 14(2):331--332, 1976.

\bibitem{Bass68}
Hyman Bass.
\newblock {\em Algebraic ${K}$-theory}.
\newblock W. A. Benjamin, Inc., New York-Amsterdam, 1968.

\bibitem{Bass79}
Hyman Bass.
\newblock Traces and {E}uler characteristics.
\newblock In {\em Homological group theory (Proc. Sympos., Durham, 1977)},
  pages 1--26. Cambridge Univ. Press, Cambridge, 1979.

\bibitem{BerrickChatterjiMislin}
A.~J. Berrick, Chatterji I., and Mislin G.
\newblock From acyclic groups to the {B}ass conjecture for amenable groups.
\newblock preprint.

\bibitem{Cliff80}
Gerald~H. Cliff.
\newblock Zero divisors and idempotents in group rings.
\newblock {\em Canad. J. Math.}, 32(3):596--602, 1980.

\bibitem{Dornhoff72}
Larry Dornhoff.
\newblock {\em Group representation theory. {P}art {B}: {M}odular
  representation theory}.
\newblock Marcel Dekker Inc., New York, 1972.
\newblock Pure and Applied Mathematics, 7.

\bibitem{Eckmann86}
Beno Eckmann.
\newblock Cyclic homology of groups and the {B}ass conjecture.
\newblock {\em Comment. Math. Helv.}, 61(2):193--202, 1986.

\bibitem{EpsteinShub68}
David Epstein and Michael Shub.
\newblock Expanding endomorphisms of flat manifolds.
\newblock {\em Topology}, 7:139--141, 1968.

\bibitem{Farkas75}
Daniel~R. Farkas.
\newblock Miscellany on {B}ieberbach group algebras.
\newblock {\em Pacific J. Math.}, 59(2):427--435, 1975.

\bibitem{FarkasMarciniak82}
Daniel~R. Farkas and Zbigniew~S. Marciniak.
\newblock Lifting idempotents in group rings.
\newblock {\em J. Pure Appl. Algebra}, 25(1):25--32, 1982.

\bibitem{FarrellHsiang70}
F.~T. Farrell and W.-C. Hsiang.
\newblock A formula for ${K}\sb{1}{R}\sb{\alpha }\,[{T}]$.
\newblock In {\em Applications of Categorical Algebra (Proc. Sympos. Pure
  Math., Vol. XVII, New York, 1968)}, pages 192--218. Amer. Math. Soc.,
  Providence, R.I., 1970.

\bibitem{FarrellHsiang78}
F.~T. Farrell and W.~C. Hsiang.
\newblock The topological-{E}uclidean space form problem.
\newblock {\em Invent. Math.}, 45(2):181--192, 1978.

\bibitem{FarrellHsiang81}
F.~T. Farrell and W.~C. Hsiang.
\newblock The {W}hitehead group of poly-(finite or cyclic) groups.
\newblock {\em J. London Math. Soc. (2)}, 24(2):308--324, 1981.

\bibitem{FarrellJones91}
F.~T. Farrell and L.~E. Jones.
\newblock Stable pseudoisotopy spaces of compact non-positively curved
  manifolds.
\newblock {\em J. Differential Geom.}, 34(3):769--834, 1991.

\bibitem{FarrellJones93}
F.~T. Farrell and L.~E. Jones.
\newblock Isomorphism conjectures in algebraic ${K}$-theory.
\newblock {\em J. Amer. Math. Soc.}, 6(2):249--297, 1993.

\bibitem{FarrellJones95}
F.~T. Farrell and L.~E. Jones.
\newblock The lower algebraic ${K}$-theory of virtually infinite cyclic groups.
\newblock {\em $K$-Theory}, 9(1):13--30, 1995.

\bibitem{FarrellLinnell}
F.~Thomas Farrell and Peter~A. Linnell.
\newblock {$K$}-theory of solvable groups.
\newblock Proc. London. Math. Soc., to appear.

\bibitem{Gersten72}
S.~M. Gersten.
\newblock On the spectrum of algebraic ${K}$-theory.
\newblock {\em Bull. Amer. Math. Soc.}, 78:216--219, 1972.

\bibitem{Goodearl76}
K.~R. Goodearl.
\newblock {\em Ring theory}.
\newblock Marcel Dekker Inc., New York, 1976.
\newblock Nonsingular rings and modules, Pure and Applied Mathematics, No. 33.

\bibitem{GrigorchukZuk01}
Rostislav~I. Grigorchuk and Andrzej {\.Z}uk.
\newblock The lamplighter group as a group generated by a 2-state automaton,
  and its spectrum.
\newblock {\em Geom. Dedicata}, 87(1-3):209--244, 2001.

\bibitem{Jacobson64}
Nathan Jacobson.
\newblock {\em Structure of rings}.
\newblock American Mathematical Society, Providence, R.I., revised edition,
  1964.

\bibitem{KrophollerLinnell88}
P.~H. Kropholler, P.~A. Linnell, and J.~A. Moody.
\newblock Applications of a new ${K}$-theoretic theorem to soluble group rings.
\newblock {\em Proc. Amer. Math. Soc.}, 104(3):675--684, 1988.

\bibitem{Lafforgue}
V.~Lafforgue.
\newblock {$K$}-th{\'{e}}orie bivariante pour les alg{\`{e}}bres de {B}anach et
  conjecture de {B}aum-{C}onnes.
\newblock preprint.

\bibitem{Linnell83}
P.~A. Linnell.
\newblock Decomposition of augmentation ideals and relation modules.
\newblock {\em Proc. London Math. Soc. (3)}, 47(1):83--127, 1983.

\bibitem{Linnell93}
Peter~A. Linnell.
\newblock Division rings and group von {N}eumann algebras.
\newblock {\em Forum Math.}, 5(6):561--576, 1993.

\bibitem{Loday76}
Jean-Louis Loday.
\newblock ${K}$-th\'eorie alg\'ebrique et repr\'esentations de groupes.
\newblock {\em Ann. Sci. \'Ecole Norm. Sup. (4)}, 9(3):309--377, 1976.

\bibitem{Milnor71}
John Milnor.
\newblock {\em Introduction to algebraic ${K}$-theory}.
\newblock Princeton University Press, Princeton, N.J., 1971.
\newblock Annals of Mathematics Studies, No. 72.

\bibitem{Quinn79}
Frank Quinn.
\newblock Ends of maps. {I}.
\newblock {\em Ann. of Math. (2)}, 110(2):275--331, 1979.

\bibitem{Quinn82}
Frank Quinn.
\newblock Ends of maps. {I}{I}.
\newblock {\em Invent. Math.}, 68(3):353--424, 1982.

\bibitem{Quinn85}
Frank Quinn.
\newblock Algebraic ${K}$-theory of poly-(finite or cyclic) groups.
\newblock {\em Bull. Amer. Math. Soc. (N.S.)}, 12(2):221--226, 1985.

\bibitem{Quinn85x}
Frank Quinn.
\newblock Geometric algebra.
\newblock In {\em Algebraic and geometric topology (New Brunswick, N.J.,
  1983)}, pages 182--198. Springer, Berlin, 1985.

\bibitem{Robinson96}
Derek J.~S. Robinson.
\newblock {\em A course in the theory of groups}.
\newblock Springer-Verlag, New York, second edition, 1996.

\bibitem{ScottWall79}
Peter Scott and Terry Wall.
\newblock Topological methods in group theory.
\newblock In {\em Homological group theory (Proc. Sympos., Durham, 1977)},
  pages 137--203. Cambridge Univ. Press, Cambridge, 1979.

\bibitem{Swan70}
Richard~G. Swan.
\newblock {\em ${K}$-theory of finite groups and orders}.
\newblock Springer-Verlag, Berlin, 1970.
\newblock Lecture Notes in Mathematics, Vol. 149.

\bibitem{Wagoner72}
J.~B. Wagoner.
\newblock Delooping classifying spaces in algebraic ${K}$-theory.
\newblock {\em Topology}, 11:349--370, 1972.

\bibitem{Waldhausen78}
Friedhelm Waldhausen.
\newblock Algebraic ${K}$-theory of generalized free products. {I}{I}{I},
  {I}{V}.
\newblock {\em Ann. of Math. (2)}, 108(2):205--256, 1978.

\bibitem{Weibel}
Charles~A. Weibel.
\newblock An introduction to algebraic ${K}$-theory.
\newblock A forthcoming graduate textbook, see
  \verb+http://www.math.rutgers.edu/~weibel/Kbook.html+.

\end{thebibliography}

\end{document}